\documentclass[11pt]{article}
\usepackage{graphicx, color, amsfonts}
\usepackage{graphicx}

\def\ds{\displaystyle}
\def\forall{\hbox{for all}~}
\def\L{{\bf L}}
\def\D{{\cal D}}
\def\E{{\cal E}}

\def\bfw{{\bf w}}

\def\ve{\varepsilon}

\def\dint{\int\!\!\!\int}

\def\N{{\mathbb N}}
\def\R{{\mathbb R}}%{I\!\!R}

\def\implies{\Longrightarrow}
\def\vp{\varphi}

\def\bfI{{\bf I}}
\def\tv{\hbox{\rm Tot.Var.}}
\def\osc{\hbox{\rm Osc.}}
\def\vs{\vskip 2em}

\def\v{\vskip 1em}
\def\O{{\cal O}}

\def\begi{\begin{itemize}}
\def\endi{\end{itemize}}
\def\C{{\cal C}}
\def\F{{\cal F}}

\def\ov{\overline}
\def\Tilde{\widetilde}
\def\Hat{\widehat}

\def\bega{\begin{array}}
\def\enda{\end{array}}
\def\meas{\hbox{meas}}

\def\bel{\begin{equation}\label}
\def\eeq{\end{equation}}
\def\sqr#1#2{\vbox{\hrule height .#2pt
\hbox{\vrule width .#2pt height #1pt \kern #1pt
\vrule width .#2pt}\hrule height .#2pt }}
\def\square{\sqr74}
\def\endproof{\hphantom{MM}\hfill\llap{$\square$}\goodbreak}

\newtheorem{theorem}{Theorem}[section]

\newtheorem{lemma}{Lemma}[section]

\newtheorem{remark}{Remark}[section]
\newtheorem{definition}{Definition}[section]

\begin{document}
\title{\bf A posteriori Error Estimates for Numerical Solutions to
Hyperbolic Conservation Laws}

\vs
\author{Alberto Bressan, Maria Teresa Chiri, and Wen Shen\\
\,
\\
Department of Mathematics, Penn State University \\
University Park, Pa.~16802, USA.\\
\,
\\
e-mails: axb62@psu.edu, mxc6028@psu.edu, wxs27@psu.edu.
}
\maketitle

\begin{abstract} The paper is concerned with {\it a posteriori} error bounds
for a wide class of numerical schemes, for
$n\times n$ hyperbolic conservation laws in one space dimension.
These estimates are achieved by a ``post-processing algorithm",
checking that the numerical solution retains small total variation, and 
computing its oscillation on suitable subdomains.
The results apply, in particular, to solutions obtained by the Godunov or the Lax-Friedrichs scheme,
backward Euler approximations, and the method of periodic smoothing.
Some numerical implementations are presented. 
\end{abstract}

\v
\section{Introduction}
\label{s:1}
\setcounter{equation}{0}
Consider the Cauchy problem for a strictly hyperbolic system of conservation laws
in one space dimension:
\bel{1}
u_t+f(u)_x\,=\,0,\eeq
\bel{2}u(0,x)\,=\,\bar u(x).\eeq
%As customary, we assume that each characteristic field is either linearly degenerate or genuinely nonlinear.
For initial data with  small total variation,
it is well known that this problem has a unique entropy-weak 
solution, depending Lipschitz 
continuously
on the initial data $\bar u$ in the $\L^1$ norm \cite{Bbook, Btut, Daf, HR}.

A closely related question is the stability and convergence of various types
of approximate solutions.
Estimates on the convergence rate for a deterministic version of the Glimm scheme 
\cite{Glimm, Liu}  were derived in \cite{BM}, and more recently in \cite{AM, BiM} for a wider
class of flux functions.
For vanishing viscosity approximations
\bel{3}
	u_t+f(u)_x\,=\,\ve\, u_{xx}\,,
\eeq
uniform BV bounds, stability and convergence as $\ve\to 0$
were proved in \cite{BiB}, while convergence rates were later established  in 
\cite{BHWY, BY04}.
Further convergence results were  proved by Bianchini 
for approximate solutions constructed by the semidiscrete (upwind)
Godunov scheme \cite{Bi03},
 and by the  Jin-Xin relaxation model  \cite{Bi04}.

A major remaining open problem is the convergence of fully discrete approximations,
such as the Lax-Friedrichs or the Godunov scheme \cite{God, HR, LV}.  
Indeed, the convergence results known for these numerical algorithms
rely on compensated compactness~\cite{DCL}.  They apply only to $2\times 2$ 
systems, and do not yield information about uniqueness or convergence rates.

For a particular class of systems, the convergence of Godunov approximations 
was proved in \cite{BJ},
relying on uniform bounds on the total variation.  
For general hyperbolic systems, however, it is known that
the Godunov scheme is unstable w.r.t.~the BV norm.  
In \cite{BBJ} an example was constructed, showing that  the total variation of a numerical solution can become arbitrarily large as $t\to +\infty$.   
Indeed, if the exact solution contains a shock with speed close to a rational multiple
of the grid size $\Delta x/\Delta t$, this can cause resonances, producing a large
amount of downstream oscillations.

Without an {\it a priori} bound on the total variation, one cannot compare an approximate solution with trajectories of the semigroup of exact solutions, 
and all the uniqueness arguments developed
in \cite{BG, BLF, BL} break down.  The counterexample in \cite{BBJ} 
can thus be regarded as a fundamental obstruction toward the derivation of {\it a priori} error estimates
for fully discrete numerical schemes.

%the existence of a convergent subsequence can only be proved by a compensated compactness argument \cite{DCL}, restricted to $2\times 2$ systems. This does not provide any information about uniqueness or continuous dependence on initial data.

To make progress, in this paper we shift our point of view,  focusing on 
{\it a posteriori} error estimates.  Namely, we assume
that an approximate solution to (\ref{1})-(\ref{2})
has been constructed by some numerical
algorithm.
Based on some additional information about the approximate solution, 
we seek an estimate
on the difference 
\bel{diff}\|u^{approx}(T, \cdot) - u^{exact}(T, \cdot)\|_{\L^1(\R)}\,.\eeq

For any sufficiently small BV initial data $\bar u$, 
it is well known that the unique entropy-admissible BV solution of (\ref{1})-(\ref{2})
has two key properties \cite{Bbook}:
\begi
\item[(i)] The total variation of $u(t,\cdot)$ remains uniformly small, for all $t\geq 0$.
\item[(ii)] Given a threshold $\rho>0$, one can identify a finite number of 
curves in the $t$-$x$ plane (shocks or contact discontinuities) such that, outside these curves,
the solution has local oscillation $<\rho$.
\endi 
The counterexample in \cite{BBJ} shows that, for an approximation constructed by the 
Godunov scheme, the property (i) sometimes can fail.   Roughly speaking, 
the result we want to prove in the present paper is the following.
Let $u^{approx}$ be an approximate solution 
produced by a conservative scheme which dissipates entropy, and assume that:
\begi
\item[(i$'$)] The total variation of $u^{approx}(t,\cdot)$ remains small, for all $t\in [0,T]$
\item[(ii$'$)] Outside a finite number of 
narrow strips in the domain $[0,T]\times\R$,
the local oscillation of $u^{approx}$ remains small.
\endi
Then the $\L^1$ distance (\ref{diff}) is small.

We emphasize that both conditions (i$'$)-(ii$'$) refer to the output of a numerical computation.
In (ii$'$), we expect that the finitely many strips where the oscillation of 
$u^{approx}$ is large
will have the form
$$
\Big\{(t,x)\,;~~t\in [a_i, b_i],~
x\in [\gamma_i(t)-\delta,\, \gamma_i(t)+\delta]\Big\},$$
where the curve $t\mapsto \gamma_i(t)$ traces the approximate location of a large shock (or a contact discontinuity) in the exact  solution.
It is also worth noting that our estimates do not require any regularity of the 
exact solution. In particular, $u^{exact}$
may well have a dense set of discontinuities.
\v
Our goal is to prove error bounds which can be applied to a wide class of approximation schemes. For future reference, we collect the basic assumptions on the system (\ref{1}), 
and the properties of the approximate solutions that will be used.
\begi
\item[{\bf (A1)}] The system (\ref{1}) is strictly hyperbolic, with 
each characteristic field being either linearly degenerate or genuinely nonlinear.
It generates a  semigroup of entropy weak solutions
$S:\D\times [0, +\infty[\,\mapsto\D$, where $\D\subset\L^1(\R;\,\R^n)$ is a 
domain containing all functions with sufficiently small total variation, namely
\bel{TVd0}
\bar u\in \L^1(\R;\,\R^n), \quad \tv\{\bar u\}~\leq \delta_0\qquad\implies
\qquad \bar u\in \D.\eeq
There exist Lipschitz constants $C_0, L_0$ such that
\bel{lip1} 
\|S_t u -S_s u\|_{\L^1}~\leq~C_0 \cdot \tv\{u\}\cdot |t-s|,\eeq
\bel{lip2}
\|S_t u - S_t v\|_{\L^1}~\leq~L_0 \|u-v\|_{\L^1}\,,\eeq
for all $u,v\in \D$ and $0\leq s\leq t$.
\item[{\bf (A2)}] For each genuinely nonlinear field, there exists a strictly convex entropy $\eta$, 
with entropy flux $q$, which  selects the admissible shocks.
\endi
We recall that the  existence of a 
semigroup generated by (\ref{1}) was proved in \cite{BiB, BC95, BCP, BLY}, 
in various degrees of generality.   In particular, it is known that the trajectories of the semigroup are the unique limits of vanishing viscosity approximations.
 To explain the additional assumption {\bf (A2)}, let $u^-, u^+$ be any two states connected by a 
 genuinely nonlinear shock
 with speed $\lambda$, so that the Rankine-Hugoniot conditions hold:
 $$\lambda\,(u^+-u^-)~=~f(u^+)-f(u^-).$$
Then, if the shock is NOT admissible, we require that the corresponding entropy should 
be 
strictly increasing, namely
\bel{enti}\lambda\left(\eta(u^-)-\eta(u^+)\right)-\left(q(u^-)-q(u^+)\right)~>~c_0 |u^--u^+|^3, 
\eeq
 for some constant $c_0>0$.
 
In the following, we shall use
test functions $\vp\in \C_c(\R^2)$ which are Lipschitz continuous with compact support, with 
Sobolev norm
\bel{lipn}\|\vp\|_{W^{1,\infty}}~\doteq~\max\bigl\{ \|\vp\|_{\L^\infty}, \,\|\vp_t\|_{\L^\infty},\, \|\vp_x\|_{\L^\infty}\bigr\}.\eeq
Given $\ve>0$, we consider approximate solutions 
$u:[0,T]\mapsto\D$ of the system of conservation laws (\ref{1}), 
taking values inside the domain of the semigroup $S$.  We assume that these solutions are 
 inductively defined
for a discrete set of times $\tau_j=j\ve$. For  $t\in [\tau_j, \tau_{j+1}[\,$ one can  then define
$u(t,\cdot)$ to be the exact solution to (\ref{1}) which coincides with $u(\tau_j,\cdot)$ at time $t=\tau_j$.
In alternative, sometimes it is more convenient to simply define $u(t,\cdot)=u(\tau_j,\cdot)$ for $t\in [\tau_j, \tau_{j+1}[\,$.

While we do not specify any particular method to construct these approximate 
solutions, two basic properties will be 
assumed.   The first is the Lipschitz continuity of the map
$t\mapsto u(t,\cdot)\in \L^1(\R\,;\R^n)$, restricted to the discrete set of times $\tau_j=j\ve$.    
The second is an approximate 
weak form of the conservation equations and the entropy conditions.   In the following, $L,C$ denote suitable constants. Moreover, the notation
$\ve\N~\doteq~\{j\ve\,;~j= 0,1,2,\ldots\}$ will be used. 
\v
\begi
\item[{\bf (AL)}] {\it 
For every $0\leq \tau<\tau' \leq T$ with $\tau,\tau'\in \ve\N$, one has}
\bel{LLip} \|u(\tau',\cdot) - u(\tau, \cdot)\|_{\L^1}~\leq~L\,(\tau'-\tau)\cdot \sup_{t\in 
[\tau, \tau']} \tv\bigl\{ u(t,\cdot)\bigr\}.\eeq

\item[{\bf (P$_\ve$)}] {\it 
For every $0\leq \tau<\tau' \leq T$ with $\tau,\tau'\in \ve\N$, and every test function $\vp\in \C^1_c(\R^2)$, one has
\bel{q1}\bega{l}\ds \left| \int u(\tau,x)\vp(\tau,x)\, dx-\int u(\tau',x)\vp(\tau,x)\, dx+
\int_\tau^{\tau'}\! \!\int \bigl\{ u\vp_t+f(u)\vp_x\bigr\}\,dx\,dt\right|\\[4mm]
\qquad\qquad\ds \leq~C \ve \|\vp\|_{W^{1,\infty}}
\cdot(\tau'-\tau)\cdot \sup_{t\in 
[\tau, \tau']} \tv\bigl\{ u(t,\cdot)\bigr\}.\enda\eeq
Moreover, assuming $\vp\geq 0$, one has the entropy inequality
\bel{q2}
\bega{l}\ds \int \eta(u(\tau,x))\vp(\tau,x)\, dx-\int \eta(u(\tau',x))\vp(\tau',x)\, dx+
\int_\tau^{\tau'}\! \!\int \bigl\{\eta(u)\vp_t+q(u)\vp_x\bigr\}\,  dxdt
\\[4mm]\qquad\qquad\ds \geq~-C \ve \|\vp\|_{W^{1,\infty}}
\cdot(\tau'-\tau)\cdot \sup_{t\in 
[\tau, \tau']} \tv\bigl\{ u(t,\cdot)\bigr\}.\enda\eeq
 }
\endi
We remark that, for an exact solution, the left hand side of  (\ref{q1}) 
would be zero, while the left hand side of  (\ref{q2}) would be non-negative.
Since here we are dealing with $\ve$-approximate solutions, we allow an error
that decreases with $\ve$, but increases with the Lipschitz constant of the test function
$\vp$.

In the present paper, two main questions will be addressed:
\begi
\item Given an approximate solution $u$ of (\ref{1})-(\ref{2}) 
satisfying {\bf (AL)} 
and {\bf (P$_\ve$)}, can one estimate the distance between $u$ and the exact solution ?
\item What kind of approximation schemes satisfy the conditions {\bf (AL)} 
and {\bf (P$_\ve$)} ? 
\endi

To answer the first question, using 
a technique introduced in \cite{BGlimm}, two types of estimates
will be derived.
\begi\item[-] On regions where the oscillation is small, the 
approximate solution $u$ is compared with the solution to a linear hyperbolic problem with constant coefficients.
\item[-] Near a point where a large jump occurs,  $u$ is compared
with the solution to a Riemann problem.
\endi
We recall that, for exact solutions, this technique yields the identity
$u(t,\cdot)=S_t \bar u$, proving that an exact solution is unique and
coincides with the corresponding semigroup trajectory 
\cite{BGlimm, Bbook, BG, BLF, BL}. 
In Sections~\ref{s:2} to \ref{s:4} we develop similar estimates
in the case of an approximate solution $u$, where the right hand side of 
(\ref{q1})-(\ref{q2}) is not zero, but vanishes of order 
$\O(1)\cdot \ve \|\vp\|_{W^{1,\infty}}$. This will provide a bound
on the difference (\ref{diff}).

An important aspect must be mentioned here. The uniqueness proofs
in \cite{BG, BLF, BL} require some additional regularity condition, such as 
``Tame Variation" or ``Tame Oscillation".    These conditions are
always satisfied by solutions constructed by front tracking or by the Glimm scheme,
but may fail for a numerically approximated solution. 
To derive rigorous error bounds, we must check that an equivalent condition is satisfied.

For a numerically computed approximation, in Section~\ref{s:5}
we introduce a
{\em post-processing algorithm}, which accomplishes three main tasks:
\begi
\item[(1)] Check that the total variation remains bounded.
\item[(2)] Trace the location of a finite number of large shocks.
\item[(3)] Check that the oscillation of the solution remains small, 
on a finite number of polygonal domains, away from the large shocks.
\endi
Step (1) is the simplest, yet the crucial one. If the total variation becomes too large, 
at some time $t$ the approximate solution $u(t,\cdot)$ can fall outside the domain $\D$ of the semigroup.
When this happens, the algorithm stops and no error estimate is achieved.  

In the favorable case where the total variation remains small, the algorithm can then proceed with steps (2) and (3).
To implement these steps, one needs to introduce certain parameters, such as the minimum size of the shocks 
which will be traced,
and the length of the time intervals $[t_j, t_{j+1}]$ used in a new partition of $[0,T]$.   
For every choice of these parameter values, the algorithm yields an error bound.  
In practice, the accuracy of this estimate largely depends on the choice of these values.
At the end of Section~\ref{s:3}, and then again at the end of 
Section~\ref{s:5}, we discuss how to choose these parameters, and 
the expected order of magnitude of the corresponding error bounds.

To complete our program, in Section~\ref{s:6}
we consider various approximation schemes, and prove that they all satisfy the properties 
{\bf (AL)}  and {\bf (P$_\ve$)}.  In particular, 
our analysis applies to: (i) Godunov's scheme, (ii) the Lax-Friedrichs' scheme, (iii) 
backward Euler  approximations, and (iv) approximate solutions obtained by periodic mollifications.

Finally, in Section~\ref{s:7} we discuss  details of the post-processing algorithm, and
present a numerical simulation.  For the ``p-system", describing isentropic 
gas dynamics in Lagrangian coordinates, we consider initial data  generating two centered rarefactions,
and two
shocks that 
eventually cross each other.  After computing an approximate solution by the Godunov scheme,
we implement the post-processing algorithm.  The two shocks are traced (as long as they remain well separated), and the remaining domain is
covered by trapezoids where the numerical solution has small oscillation (away from interaction times).

\section{Solutions with small oscillation} \label{s:2}
\setcounter{equation}{0}
In this section we  begin by studying the case where no large shocks are present.   
Let $u=u(t,x)$ be an approximate solution which satisfies {\bf (AL)} and {\bf (P$_\ve$)}.
Consider an open interval $]a,b[\,$,
fix a point $\xi$ with $a<\xi<b$ %\in \,]a,b[\,$, 
and set 
\bel{Adef}A ~=~ Df(u(0,\xi)).\eeq
Assuming that all characteristic speeds satisfy
\bel{speed}\lambda^-\,<\,\lambda_i(u)\,<\,\lambda^+,\qquad i=1,\ldots,n,\eeq
fix $\tau\in \ve\N$ and consider the trapezoidal domain
\bel{trap}\Delta~=~\Big\{(t,x)\,;~t\in [0,\tau], ~~a(t)~\doteq~ a+\lambda^+ t ~<~x
<~b+\lambda^- t~\doteq~ b(t)\Big\}.\eeq
Following an approach introduced in \cite{BGlimm}, error estimates will be obtained
by comparing $u$   with the solution $w$
of the linear hyperbolic system with constant coefficients
\bel{le}
w_t + A w_x~=~0,\qquad\qquad w (0,x) = u (0,x).\eeq

For this purpose, let $\{\ell_1,\ldots,\ell_n\}$ and $\{ r_1,\ldots,r_n\}$
be dual bases of left and right eigenvectors of the matrix $A$, normalized so that
\bel{rli}|r_i|~=~1,\qquad\qquad \ell_i\cdot r_j~=~\delta_{ij}~=~\left\{\bega{rl}
1\quad\hbox{if}~~i=j,\cr
0\quad\hbox{if}~~i\not= j.\enda\right.\eeq
Let $\lambda_1,\ldots,\lambda_n$ be the corresponding eigenvalues of $A$.
For each $i$, consider the scalar functions
$$u_i(t,x)~=~\ell_i\cdot u(t,x),\qquad\qquad w_i(t,x)~=~\ell_i\cdot w(t,x).$$
By (\ref{le}), $w_i$ solves the scalar  linear  equation
$$w_{i,t}+\lambda_i w_{i,x}~=~0,\qquad\qquad w_i(0,x)~=~u_i(0,x).$$
%In the end, this amounts to estimate by how much each component 
%$u_i$ fails to satisfy the scalar conservation law
%\bel{ui1}u_i + (\lambda_i u_i)_x~=~0.\eeq
For each $i=1,\ldots,n$,
we will estimate  the difference
$u_i(\tau,\cdot)-w_i(\tau,\cdot)$.
\v
As a preliminary, consider a BV function $g:[\alpha,\beta]\mapsto \R$.
Since $g$ is regulated, it admits left and right limits $g(x-)$, $g(x+)$
at every point $x$.
By possibly modifying $g$ on the  countable set where it has jumps,
 we can assume that
\bel{gj} g(x+)\cdot g(x-)\,\leq \,0\qquad\implies\qquad  g(x)\,=\,0. \eeq
We can then select countably many maximal open subintervals 
$]a_j, b_j[\, \subset [\alpha,\beta]$ where $g$ has constant sign.
Namely,
\begi
\item[{\bf (G)}] {\it 
 $g$ has constant sign on each $\,]a_j,b_j[\,$, and changes sign on 
 every  neighborhood of
each endpoint
$a_j, b_j$ (unless $a_j=\alpha$ or $b_j=\beta$).  
Moreover, $g(x)=0$ for $x\notin \bigcup_j [a_j, b_j]$.}
\endi
For a given $\ve>0$, 
consider the test function with Lipschitz constant $\|\phi_x\|_{\L^\infty} = \ve^{-2/3}$
\bel{phidef}\phi(x)~\doteq~\left\{\bega{cl} \ds\min\left\{1,~ {x-a_j\over \ve^{2/3}}, ~{b_j-x\over\ve^{2/3}} \right\}
\quad &\hbox{if $x\in [a_j,b_j]$ and  $g$ is positive on $]a_j, b_j[\,$,}\\[4mm]\ds
\max\left\{-1,~ {a_j-x\over \ve^{2/3}},~ {x-b_j\over\ve^{2/3}} \right\}
\quad &\hbox{if $x\in [a_j,b_j]$ and  $g$ is negative on $]a_j, b_j[\,$,}\\[4mm]
0\quad &\hbox{if $x\notin \cup_j [a_j, b_j]$.}
\enda\right.
\eeq

\begin{figure}[ht]
\centerline{\hbox{\includegraphics[width=10cm]{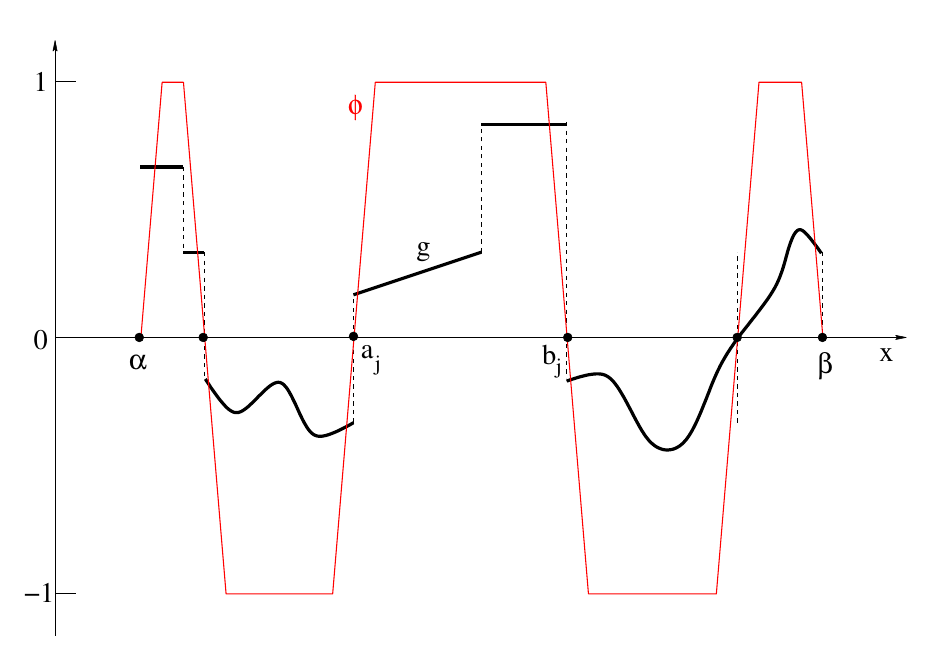}}}
\caption{\small  The test function $\phi$
defined at (\ref{phidef}), with Lipschitz constant $\ve^{-2/3}$.
}
\label{f:hyp212}
\end{figure}

\begin{lemma}\label{l:21}
Let $g:[\alpha, \beta]\mapsto\R$ be as above.
If $g$ is strictly positive (or strictly negative) for all $x\in [\alpha,\beta]$,
then  
\bel{ge1}
\int_{\alpha+\ve^{2/3}}^{\beta-\ve^{2/3}} |g(x)|\, dx~\leq~\int_\alpha^\beta \phi(x)\, g(x)\, dx.\eeq
 On the other hand, if $g(\xi)=0$ for some $\xi\in [\alpha,\beta]$, then
 \bel{ge2}
\int_\alpha^\beta |g(x)|\, dx~\leq~\int_\alpha^\beta \phi(x)\, g(x)\, dx +
 2\ve^{2/3}\cdot  \tv\bigl\{ g\,;~[\alpha,\beta]\bigr\} .
\eeq
\end{lemma}
\v
{\bf Proof.} {\bf 1.} If $g$ has always the same sign, then by construction
$\phi(x)\,g(x)\geq 0$ for all $x$, while
$\phi(x)=\hbox{sign}(g(x))$ for $x\in [\alpha+\ve^{2/3},\, \beta-\ve^{2/3}]$.
Hence the estimate (\ref{ge1}) is trivially true.
%$$\int_\alpha^\beta |g(x)|\, dx-\int_\alpha^\beta \phi(x)\, g(x)\, dx
%~=~\int_\alpha^{\alpha+ \ve^{1/3}} {x-\alpha\over\ve^{1/3}} |g(x)|\, dx +
%\int_{\beta- \ve^{1/3}}^\beta {\beta-x\over\ve^{1/3}} |g(x)|\, dx ~\leq~\ve^{1/3}\cdot
%\|g\|_{\L^\infty}\,.$$
\v
{\bf 2.} If $g$ changes sign, consider the maximal subintervals 
$[a_j, b_j]$ where $g$ has a constant sign, as in {\bf (G)}. We then have the estimate
$$\bega{l} \ds\int_\alpha^\beta |g(x)|\, dx~=~\sum_j \int_{a_j}^{b_j} |g(x)|\, dx\\[4mm]
\ds\quad\leq~\sum_{b_j-a_j\leq 2\ve^{2/3}} (b_j-a_j) \sup_{a_j<x<b_j} |g(x)|
\\[4mm]
\ds\qquad\qquad\qquad
+ \sum_{b_j-a_j\geq 2 \ve^{2/3}}\left(\int_{a_j+\ve^{2/3}}^{b_j-\ve^{2/3}} +
\int_{a_j}^{a_j+\ve^{2/3}} + \int_{b_j-\ve^{2/3}}^{b_j}
\right) |g(x)|\, dx
\\[4mm]
\ds\quad\leq~2\ve^{2/3}\cdot\!\!\!\!\! \sum_{b_j-a_j\leq 2\ve^{2/3}} \tv\{ g\,;~]a_j, b_j]\}
%\\[4mm]\ds\qquad\qquad \qquad  
+\int_\alpha^\beta \phi(x)\, g(x) \, dx + 2\ve^{2/3}\cdot\!\!\!\!\! \sum_{b_j-a_j> 2\ve^{2/3}}
\sup_{a_j<x<b_j} |g(x)|\\[4mm]
\ds\quad =~
\int_\alpha^\beta \phi(x)\, g(x) \, dx +2\ve^{2/3}\cdot  \tv\bigl\{ g\,;~[\alpha,\beta]\bigr\} .
\enda
$$
\endproof

\begin{remark}
{\rm Here and in the sequel, one could prove similar results by replacing the exponent $2/3$ with any number $\gamma\in \,]0, 1[\,$, and working with 
test functions which are Lipschitz continuous with constant $\ve^{-\gamma}$.   
Our choice of $\gamma=2/3$ is motivated by the heuristic expectation that, 
in most cases, this should yield the sharpest error bounds.   
See Remark~\ref{r:31} for further discussion of this point.}
\end{remark}

We can now state the main result of this section, providing an error estimate on the trapezoidal domain (\ref{trap}).

\begin{lemma}\label{l:22} There exists a constant $C_1$ such that the following holds.  For a given $\ve>0$,
let $u$ be an approximate solution of (\ref{1})  that satisfies the 
%properties {\bf (AL)} and 
property {\bf (P$_\ve$)}. Let $\Delta$ be the trapezoid in (\ref{trap}), and let $w$ be the solution to the linear Cauchy problem
(\ref{le}), with $A$ as in (\ref{Adef}). Then
% \small 
\bel{ee3}\bega{l}
\ds
\int_{a+\tau \lambda^+ +\ve^{2/3}}^{b+\tau \lambda^--\ve^{2/3} } 
\bigl| u(\tau,x) - w(\tau,x)\bigr|\, dx\\[4mm]
\ds \quad  \leq~C_1\bigg(\tau\cdot \sup_{(t,x)\in \Delta}\bigl|u(t,x)-u(0,\xi)|
+\tau \ve^{1/3}+\ve^{2/3}\bigg)\cdot \sup_{t\in 
[0,\tau]} \tv\bigl\{ u(t,\cdot)\bigr\}.\enda
\eeq
\end{lemma}

%\begin{figure}[ht]
%\centerline{\hbox{\includegraphics[width=13cm]{FIG/hyp211.pdf}}}
%\caption{\small  The construction used in the proof of Lemma~\ref{l:22}.
%}
%\label{f:hyp211}
%\end{figure}

{\bf Proof.} {\bf 1.} Fix $i\in \{1,\ldots,n\}$. On the interval $[\alpha,\beta]
\doteq [a+\tau\lambda^+,\, b+\tau \lambda^-]$, consider the scalar function
\bel{gdef}g_i(x) ~\doteq~\ell_i \cdot \bigl[ u(\tau, x) - w(\tau,x)\bigr]~=~
\ell_i \cdot \bigl[ u(\tau, x) - u(0,x-\lambda_i\tau)\bigr].\eeq

Let  $\phi_i:[\alpha,\beta]\mapsto [-1,1]$ be 
the  function with Lipschitz constant $\|\phi_i'\|_{\L^\infty}=\ve^{-2/3}$,
defined  as in (\ref{phidef}) with $g$ replaced by $g_i$.  We then extend 
$\phi_i$ to the entire real line by setting $\phi_i(x)=0$ if $x\notin [\alpha,\beta]$, 
and consider a  test function $\vp_i=\vp_i(t,x)$ such that
\bel{tnu}\vp_i(t,x)~=~\phi_i\bigl(x- \lambda_i(t-\tau)\bigr)\,\ell_i\qquad
\quad\hbox{for} \quad t\in [0,\tau], ~x\in\R\,.\eeq
%\left\{\bega{cl} \min\Big\{ 1,\, \nu x,\, \nu(\tau-x)\Big\}\cdot \phi\bigl(x- \lambda_i(t-\tau)\bigr)\,\ell_i
%\qquad &\hbox{if} \quad t\in  [0,\tau],\\[3mm]
%0\qquad &\hbox{if}\quad t\notin [0,\tau].
%\enda\right.\eeq
\v
{\bf 2.} Observing that  $\lambda_i\in [0,1]$ and $\|\vp_i\|_{W^{1,\infty}}=|\ell_i|\, \ve^{-2/3} $,
by (\ref{q1}) we now obtain
\bel{v2}\bega{l}\ds\left|\int \vp_i(0,x) u(0,x)\, dx - \int \vp_i(\tau,x) u(\tau,x)\, dx +\int_0^\tau
\int \bigl\{u\vp_{i,t} + f(u)\vp_{i,x}\bigr\}\, dx dt\right|\\[4mm]
\qquad\qquad\ds \leq~C\ve \|\vp_i\|_{W^{1,\infty}}\cdot \tau \cdot \sup_{t\in 
[0,\tau]} \tv\bigl\{ u(t,\cdot)\bigr\}\\[4mm]
\qquad\qquad\ds\leq~C' \ve^{1/3}\,\tau \cdot \sup_{t\in 
[0,\tau]} \tv\bigl\{ u(t,\cdot)\bigr\}.\enda\eeq
\v
{\bf 3.} For future use we observe that, if $x\mapsto u(x)$ is Lipschitz and 
$u^*=u(\xi)$ for some $\xi\in [x_1, x_2]$, then
$$\bega{rl}\ds
\int_{x_1}^{x_2} \Big|\ell_i \bigl( f(u)_x -\lambda_i  u_x\bigr)\Big|\, dx
&\ds=~\int_{x_1}^{x_2} \Big|\ell_i  \bigl[Df(u) - Df(u^*)\bigr]  u_x\Big|\, dx\\[4mm]
&\ds\leq~C_0\,\sup_{x_1<x<x_2} |u(x)-u^*|\cdot \int_{x_1}^{x_2}  |u_x|\, dx.
\enda
$$
Here $C_0$ is a constant depending only on the function $f$. 
By an approximation argument, for any BV function $x\mapsto u(x)$ we conclude
\bel{tv11}\tv\Big\{\ell_i\bigl( f(u)-\lambda_i u\bigr)\,;~[x_1,x_2]\Big\}~\leq~C_0
\left(\sup_{x_1<x<x_2} |u(x)-u^*|\right)\cdot \tv\bigl\{u\,;~[x_1, x_2]\bigr\}.
\eeq

\v
{\bf 4.} 
Since $w$ is a solution to the linear equation (\ref{le}), the choice of the test function $\vp_i$ in (\ref{tnu}) implies
\bel{v5}\int \vp_i(0,x) \,u(0,x)\, dx~=~\int \vp_i(0,x) \,w(0,x)\, dx~=~\int \vp_i(\tau,x) 
\,w(\tau,x)\, dx.
\eeq
Moreover, calling $u^*=u(0,\xi)$,
%since every function $x\mapsto \vp(t,x)$  in (\ref{test}) has compact support, 
integrating by parts and using  (\ref{tv11}) together with the bound
$\phi_i(x)\in [-1,1]$, we obtain
%and norm $\|\vp\|_{\L^\infty}\leq |\ell_i|$, we can write
\bel{v6}
\bega{l}\ds\left|\int_0^\tau
\int \bigl\{u\vp_{i,t} + f(u)\vp_{i,x}\bigr\}\, dx dt\right|~ =~\left|\int_0^\tau\int \bigl[ f( u) -\lambda_i u
\bigr]\, \vp_{i,x} \,dx dt\right|\\[4mm]
\ds\qquad \leq~\int_0^\tau \tv\Big\{ \ell_i (f( u) -\lambda_i u) \,;~[a+\lambda^+\tau+
(t-\tau)\lambda_i\,,~b+\lambda^-\tau+
(t-\tau)\lambda_i]\Big\}\, \|\vp_i\|_{\L^\infty}\,dt
\\[4mm]
\ds\qquad \leq ~
C_0\,\sup_{ (t,x)\in \Delta} 
|u(t,x)-u^*|\cdot
\int_0^\tau\tv\Big\{u(t,\cdot)\,;~\,]a+t\lambda^+,\, b+t\lambda^-[\Big\}\, dt.
\enda
\eeq
\v
{\bf 5.} By (\ref{v5}), combining (\ref{v2}) with (\ref{v6})  we conclude
\bel{v8}\bega{l}\ds~\int_{a+\tau \lambda^+}^{b+\tau\lambda^-} \phi_i(x)\,\ell_i \bigl[ w(\tau,x)-u(\tau,x)\bigr]\, dx~=~
\int \vp_i(0,x) u(0,x)\, dx - \int \vp_i(\tau,x) u(\tau,x)\, dx\\[4mm]
\qquad\ds  \leq~\left|\int_0^\tau
\int \bigl\{u\vp_{i,t} + f(u)\vp_{i,x}\bigr\}\, dx dt\right| + \E_i\,,\enda\eeq
where
\bel{Edef}\bega{rl} \E_i&\doteq~\ds
\left|\int \vp_i(0,x) u(0,x)\, dx - \int \vp_i(\tau,x) u(\tau,x)\, dx +\int_0^\tau
\int \bigl\{u\vp_{i,t} + f(u)\vp_{i,x}\bigr\}\, dx dt\right|\\[4mm]
&\ds \leq~C\ve\|\vp\|_{W^{1,\infty}} \cdot \tau\cdot \sup_{t\in [0,\tau]} 
\tv\bigl\{u(t,\cdot)\bigr\}.\enda\eeq
Notice that the above inequality follows from (\ref{q1}).  In addition, the 
 first term on the right hand side of (\ref{v8}) is estimated by (\ref{v6}).
 \v
{\bf 6.}  If the function $g_i(x)\doteq \ell_i \bigl[ w(\tau,x)-u(\tau,x)\bigr]$ 
always keeps the same sign, we now use
(\ref{ge1}). If it changes sign at least once, we use (\ref{ge2}).
Combining the two cases, by (\ref{v6}) and (\ref{v8}) we deduce
\bel{v9}\bega{l}\ds
\int_{a+\tau \lambda^++\ve^{2/3}}^{b+\tau\lambda^--\ve^{2/3} } \Big|\ell_i \bigl[ w(\tau,x)-u(\tau,x)\bigr]
\Big|\, dx
\\[4mm]
\ds\qquad \leq~C_0\,\sup_{ (t,x)\in \Delta} 
|u(t,x)-u^*|\cdot
\int_0^\tau\tv\Big\{u(t,\cdot)\,;~\,]a+t\lambda^+,\, b+t\lambda^-[\Big\}\, dt + \E_i
\\[4mm]
\qquad\qquad \ds + 2\ve^{2/3}\cdot \tv
\Big\{g_i\,;~]a+\tau\lambda^+, b+\tau\lambda^-[\Big\}.\enda
\eeq
\v
{\bf 7.}
Recalling (\ref{rli}),  for any vector $v= \sum_i c_i r_i \in\R^n$ one has
$$|v|~\leq~\sum_{i=1}^n| c_i|~=~\sum_{i=1}^n |\ell_i \cdot v|.$$
We use this inequality with $v= w(\tau,x) - u(\tau,x)$. 
Using (\ref{gdef}) to compute the total variation of $g_i$,     summing the inequalities  
(\ref{Edef})-(\ref{v9}) for $i=1,\ldots, n$,
we obtain
\bel{v11}\bega{l}\ds
\int_{a+\tau \lambda^++\ve^{1/3}}^{b+\tau\lambda^--\ve^{1/3}} \Big| w(\tau,x)-u(\tau,x)
\Big|\, dx\\[4mm]
\ds\quad \leq ~n  C_0\,\sup_{ (t,x)\in \Delta} 
|u(t,x)-u^*|\cdot \tau \cdot \sup_{t\in [0,\tau]} 
\tv\Big\{u(t,\cdot)\,;~\,]a+t\lambda^+,\, b+t\lambda^-[\Big\}\\[4mm]
\ds\qquad +nC\ve\|\vp\|_{W^{1,\infty}} \cdot \tau\cdot \sup_{t\in [0,\tau]} 
\tv\bigl\{u(t,\cdot)\bigr\} \\[4mm]
\ds\qquad +2\ve^{2/3} \left(\sum_i |\ell_i|\right)
\cdot \Big(
\tv\bigl\{ u(0,\cdot)\,;~[a,b]\bigr\} + \tv\bigl\{ u(\tau,\cdot)\,;~[a+\tau \lambda^+, ~b+\tau\lambda^-]
\bigr\}\Big).
\enda
\eeq
This yields (\ref{ee3}), for a suitable constant $C_1$.
\endproof
\section{Error bounds for solutions without large shocks}
 \label{s:3}
\setcounter{equation}{0}
%In this section we show how the estimates proved in Lemma~\ref{l:22}
%can be put to use, to derive error bounds for approximate solutions
%without large shocks.
Consider an approximate solution $u=u(t,x)$ of (\ref{1})-(\ref{2}), constructed by a numerical algorithm with time step $\ve>0$, which satisfies the properties {\bf (AL)} and {\bf (P$_\ve$)}.
%To illustrate our approach,
%in this section we shall still assume that  $u$ does not contain large shocks.  
We fix a new time step $h>\!>\ve$, and split the interval $[0,T]$ into subintervals 
$[t_j, t_{j+1}]$ with $t_j = j\, h$.  Throughout the following we choose $h\approx \ve^{1/3}$,  say
\bel{hO}
c_0\ve^{1/3}~\leq~h~\leq~\ve^{1/3}\eeq
for some constant $c_0>0$,
and assume that both $h$ and $T$ are integer multiples of $\ve$. To simplify the discussion, 
we also assume that $T=\nu h$ for some integer $\nu$.  Notice that, in the general case, 
one can consider the time $T'$ such that
$$T'~=~\nu h~\leq~ T~< ~(\nu+1) h$$ for some integer $\nu$. By (\ref{LLip}), 
the difference can then be estimated by
$$\|u(T,\cdot)- u(T',\cdot)\|_{\L^1}~\leq~L h \cdot \sup_{t\in [0,T]}\tv\{u(t,\cdot)\}~=~\O(1)\cdot \ve^{1/3}.$$

\begin{figure}[ht]
\centerline{\hbox{\includegraphics[width=12cm]{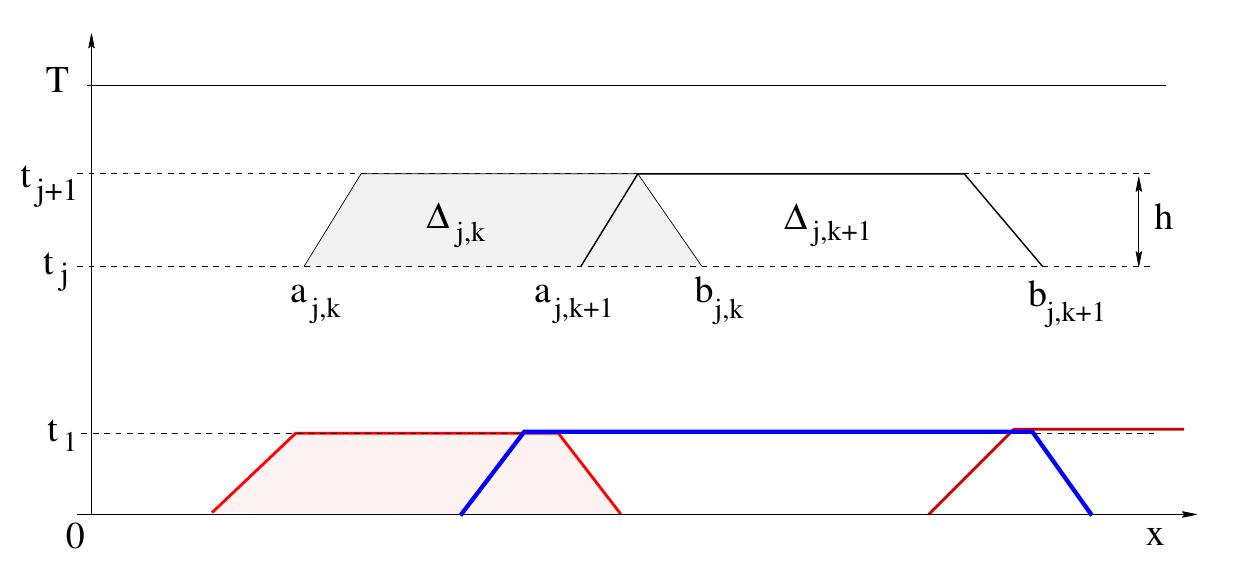}}}
\caption{\small Covering the strip $[0,T]\times \R$ with finitely many trapezoids $\Delta_{jk}$.
}
\label{f:hyp215}
\end{figure}

As shown in Fig.~\ref{f:hyp215}, for any given $j=1,2,\ldots, \nu-1$, we cover the real line with finitely many
intervals $]a_{jk}, b_{jk}[\,$, $k=1,\ldots, N(j)$, so that
$$-\infty~=~a_{j,1}~<~ a_{j,2} ~< ~b_{j,1}~<~ a_{j,3}~ <~\cdots 
~<~ a_{j,N(j)} ~< ~b_{j,N(j)-1}~<~ b_{j,N(j)}~ =~ +\infty.$$ 
We then cover each strip
$[t_j, t_{j+1}]\times \R$ with the 
trapezoids $\Delta_{jk}$, $k=1,\ldots, N(j)$. 
For convenience, these will be expressed as the convex closure of their four vertices:
$$
\Delta_{jk}~=~\hbox{co} \Big\{ (t_j ,a_{jk}), ~~(t_j, b_{jk}), ~~(t_{j+1}, a_{jk}+h\lambda^++\ve^{2/3}),~~   (t_{j+1}, b_{jk}+h\lambda^--\ve^{2/3})\Big\}.$$
Equivalently:
\bel{Djk}\bega{l}
\Delta_{jk}~=~\bigg\{ (t,x)\,;~t\in [t_j, t_{j+1}],\quad \\[4mm]
\ds\quad {t_{j+1}-t\over h} a_{jk} + {t-t_j\over h} 
(a_{jk} + h \lambda^+ + \ve^{2/3}) ~\leq~x~\leq~{t_{j+1}-t\over h} b_{jk} + {t-t_j\over h} 
(b_{jk} + h \lambda^- - \ve^{2/3}) \bigg\}.\enda\eeq

%Of course, when $k=1$ or $k=N(j)$, the sets $\Delta_{jk}$ are unbounded and need a somewhat different definition. Namely,
%$$\Delta_{j1}~=~\bigcup_{m\geq 1} \hbox{co} \Big\{ (t_j ,b_{jk}-m), ~~(t_j, b_{jk}), ~~(t_{j+1}, b_{jk}+h\lambda^--m),~~   (t_{j+1}, b_{jk}+h\lambda^- - \ve^{2/3})\Big\},$$
%while $\Delta_{j, N(j)}$ is defined similarly.

By suitably choosing the points $a_{jk}$, $b_{jk}$, we can assume that the intervals 
$$J_{jk} '~=~[a_{jk}+h\lambda^++\ve^{2/3} \,,~ b_{jk}+ h \lambda^--\ve^{2/3}], \qquad k=1,\ldots,N(j),$$
form a partition of $\R$. Namely
\bel{GG}b_{jk}+ h \lambda^--\ve^{2/3}~=~a_{j,k+1}+h\lambda^++\ve^{2/3},
\qquad\quad k=1,\ldots, N(j)-1.\eeq
Furthermore, by choosing the bases of all trapezoids to have length
\bel{lk}
b_{jk}- a_{jk}~>~2h(\lambda^+-\lambda^-),\eeq
we can assume that each point $(t,x)\in [t_j, t_{j+1}]\times\R$ is contained in at least one and in not more than two of these trapezoids.

Next, we recall that the oscillation of $u$ over a set $\Delta$ is defined as
$$\osc\bigl\{ u\,;~\Delta\bigr\}~\doteq~\sup_{(t,x), (s,y)\in \Delta}~|u(t,x)-u(t,y)|.$$
For each fixed $j\in \{1,\ldots,\nu\}$, the maximum oscillation of $u$ over all
trapezoids $\Delta_{jk}$ will be denoted by
%By choosing the domains $\Delta_{jk}$ small enough, we assume that
%\bel{tv3}\tv\bigl\{ u(t_j,\cdot)\,;~]a_{jk}, b_{jk}[\,\bigr\}~\leq~\rho.\eeq
\bel{tv4}\kappa_j~\doteq~\max_{1\leq k\leq N(j)}
\osc\bigl\{ u\,;~\Delta_{jk}\bigr\}\,.\eeq

% Applying the estimate (\ref{ee3}) to each interval
%$J_k'$, we obtain a bound on the difference   between $u(t+h, \cdot)$ and 
%an approximate solution $w(t+h,\cdot)$ obtained by patching together a family of
%solutions $w^{(i)}$ to linear problems with constant coefficients.  
%In turn, by the analysis in \cite{BGlimm}, this will yield an estimate on the difference between $u(t+h,\cdot)$ and the exact solution.

Let now $S:[0,+\infty[\,\times \D\mapsto\D$ be the Lipschitz semigroup generated by the hyperbolic system (\ref{1}), as in (\ref{lip1})-(\ref{lip2}).
In particular, $t\mapsto S_t\bar u$ yields
 the exact solution  to the Cauchy problem  (\ref{1})-(\ref{2}).   As proved in
 \cite{BGlimm, Bbook}, for any approximate solution $u$ one has 
 the error estimate
\bel{es5}
\|u(T,\cdot) - S_T\bar u\|_{\L^1}~\leq~L_0\cdot \sum_{j=0}^{\nu-1}
\Big\| u(t_{j+1}, \cdot) -S_{h} u(t_j, \cdot)\Big\|_{\L^1}\,.\eeq
For each $j$, we will show that the corresponding term on the 
right hand side of (\ref{es5})
can be estimated using (\ref{ee3}).

Consider the covering of the strip $[t_j, t_{j+1}]\times \R$
in terms of the  trapezoids $\Delta_{jk}$, introduced at  (\ref{Djk}).  
As in (\ref{le}), for $k=1,\ldots,N(j)$ we shall denote by $w^{(k)}$ the solution to the
linearized problem with constant coefficients 
\bel{lei}
w_t + A w_x~=~0,\qquad\quad w (t_j,\cdot ) = u (t_j,\cdot),
\qquad \quad A\doteq Df\bigl(u(t_j, \xi_k)\bigr),\eeq
for some given points $\xi_k\in \,]a_{jk},\,b_{jk}\,[$.

Let $\ell_i^{(k)}$ be the $i$-th left eigenvector of the above matrix $A$, normalized as in (\ref{rli}).
Using (\ref{v9}) on each trapezoid $\Delta_{jk}$ we obtain
\bel{v99}\bega{l}\ds
\int_{a_{jk}+h_j \lambda^+ +\ve^{2/3}}^{b_{jk}+h_j \lambda^--\ve^{2/3} } 
\left| \ell_i^{(k)}\cdot \bigl[u(t_{j+1},x) - w^{(k)}(t_{j+1},x)\bigr] \right|\, dx
\\[4mm]
\ds \leq~ C_0\cdot\osc\{u\,;~\Delta_{jk}\} \cdot   \int_{t_j}^{t_{j+1}}
\tv\Big\{u(t,\cdot)\,;~]a_{jk}+(t-t_j)\lambda^+, b_{jk}+(t-t_j)\lambda^-[\Big\}\, dt\\[4mm]
\quad\ds +\left|\int \vp_i^{(k)}(t_j,x) u(t_j,x)\, dx - 
\int \vp_i^{(k)}(t_{j+1},x) u(t_{j+1},x)\, dx + \dint_{\Delta_{jk}} \left\{ u\vp^{(k)}_{i,t} + f(u)\vp^{(k)}_{i,x}
\right\}\, dx dt\right|\\[4mm]
\quad+ 2\ve^{2/3}\cdot\tv\Big\{\ell_i\cdot  u(t_j,\cdot)\,;~[a_{jk},b_{jk}]\Big\} \\[4mm]
\quad+ 2\ve^{2/3} \cdot\tv\Big\{ \ell_i \cdot u(t_{j+1},\cdot)\,;~[a_{jk}+ h_j\lambda^+\,,~b_{jk}+h_j\lambda^-]
\Big\}\bigg)\\[4mm]
\doteq~A_{ik} + B_{ik} + C_{ik} + D_{ik}\,.
\enda
\eeq

For notational convenience, call $\chi_{jk}$ the characteristic function of 
the interval $[a_{jk}+h_j \lambda^+ +\ve^{2/3}\,,~b_{jk}+h_j \lambda^--\ve^{2/3}]$.
Our next goal is to estimate the quantity
\bel{Ej}
E_j~\doteq~\int_{-\infty}^{+\infty} \left| u(t_{j+1},x) - \sum_k w^{(k)}(t_{j+1},x)\cdot\chi_{jk}(x)
\right|\, dx\,.\eeq
This can of course be achieved by summing the terms 
on the right hand side of (\ref{v99}) over all $i=1,\ldots,n$
and $k=1,\ldots, N(j)$.
Toward this goal, we recall the key assumption that every point 
$(t,x)\in [t_j, t_{j+1}]\times\R$ belongs to one and no more than two of the trapezoids $\Delta_{jk}$. More precisely, we have the implication
\bel{djkk} |k-k'|\,\geq\, 2\qquad\implies\qquad 
\Delta_{jk}\cap\Delta_{jk'}~=~\emptyset.\eeq
Recalling (\ref{tv4}), for a fixed $i$ we thus obtain
\bel{Ak}
\sum_{k=1}^{N(j)} A_{ik}~\leq~C_0\, \kappa_j\cdot \int_{t_j}^{t_{j+1}} 2\,\tv\{u(t,\cdot)\}\, dt.\eeq
\bel{CDk}
\sum_{k=1}^{N(j)} C_{ik}~\leq~4\ve^{2/3}\,\tv\{u(t_j,\cdot)\},\qquad \qquad 
 \sum_{k=1}^{N(j)} D_{ik}~\leq~4\ve^{2/3}\,\tv\{u(t_{j+1},\cdot)\}.\eeq
The estimate for $\sum_k B_{jk}$ is a bit more delicate, because if we use 
(\ref{q1}) separately on each subdomain $\Delta_{jk}$, the error term
on the right side would be multiplied by $N(j)$, which can be a very large number.

For this reason, we argue as follows.  For each $i\in\{1,\ldots,n\}$, we
consider test functions 
$\vp$, $\Tilde\vp_i$ which satisfy, for  $t\in [t_j, t_{j+1}]$,
$$\vp_i(t,x)~=~\left\{\bega{cl} \vp^{(k)}_i(t,x)\qquad &\hbox{if}~~(t,x)\in \Delta_{jk}\,,~~k~\hbox{even},\cr
0 \qquad &\hbox{otherwise.}\enda\right.$$
$$\Tilde \vp_i(t,x)~=~\left\{\bega{cl} \vp^{(k)}_i(t,x)\qquad &\hbox{if}~~(t,x)\in \Delta_{jk}\,,~~k~\hbox{odd},\cr
0 \qquad &\hbox{otherwise.}\enda\right.$$
For convenience, we denote by $\ell_{\rm max}$ an upper bound for the norm of all 
left eigenvectors $\ell_i=\ell_i(u)$ of all matrices $A(u)=Df(u)$, normalized as
in   (\ref{rli}). With this notation we have
\bel{vpp}
\|\vp_i\|_{W^{1\infty}}~\leq~\ell_{\rm max}\cdot \ve^{-2/3} \,, \qquad\qquad \|\Tilde \vp_i\|_{W^{1\infty}}~\leq~\ell_{\rm max}\cdot \ve^{-2/3} .\eeq

Applying (\ref{q1}) to the test function $\vp_i$, then to $\Tilde\vp_i$, we obtain
$$\sum_{k~{\rm even}}~B_{ik}~\leq~C\ve h \,\ve^{-2/3}\ell_{\rm max} \cdot \sup_{t\in [t_{j}, t_{j+1}]}
\tv\bigl\{ u(t,\cdot)\big\},$$
$$\sum_{k~{\rm odd}}~B_{ik}~\leq~C\ve h\, \ve^{-2/3}\ell_{\rm max} \cdot \sup_{t\in [t_{j}, t_{j+1}]}
\tv\bigl\{ u(t,\cdot)\big\}.$$
Summing over $k$, we thus obtain
\bel{Bk}\sum_{k=1}^{N(j)} B_{ik}~\leq~2C h \ve^{1/3}\ell_{\rm max} \cdot \sup_{t\in [t_{j}, t_{j+1}]}
\tv\bigl\{ u(t,\cdot)\big\}.\eeq
All together, the inequalities (\ref{Ak}), (\ref{CDk}), and (\ref{Bk}) yield
\bel{ee8}
\bega{l}\ds
\int_{-\infty}^{+\infty} \left| u(t_{j+1},x) - \sum_k w^{(k)}(t_{j+1},x)\cdot\chi_{jk}(x)
\right|\, dx\\[4mm]
\qquad\ds\leq~C_0\, 2n\kappa_j\int_{t_j}^{t_{j+1}}\tv\bigl\{ u(t,\cdot)\bigr\}\, dt 
+ 4n\ve^{2/3}\tv\bigl\{ u(t_j,\cdot)\bigr\} \\[4mm]
\qquad\qquad +\ds
  4n\ve^{2/3}\tv\bigl\{ u(t_{j+1},\cdot)\bigr\} + 2C h\, \ve^{1/3}\ell_{\rm max} \cdot \sup_{t\in [t_{j}, t_{j+1}]}
\tv\bigl\{ u(t,\cdot)\big\}.
\enda
\eeq

Next, we replace the approximate solution $u$ with 
the exact solution $u^{exact}(t_j+s,\cdot)~=~S_s u(t_j,\cdot)$
of (\ref{1}) having the same data at $t=t_j$.  As proved in \cite{BGlimm},
with the same notation used in (\ref{Ej}), as long as $u(t_j,\cdot)\in\D$
remains in the domain of the semigroup, one has
\bel{ee9}\bega{l}
\ds
\int_{-\infty}^{+\infty} \left| u^{exact}(t_{j+1},x) - \sum_k w^{(k)}(t_{j+1},x)\cdot\chi_{jk}(x)
\right|\, dx\\[4mm]
%\qquad\ds 
%\leq~C_2\, \left(\max_{1\leq k\leq N(j)} \osc\{ u\,; \Delta_{jk}\}
%\right)\cdot h_j\,
%\tv\bigl\{ u(t_j,\cdot)\bigr\},
%\\[4mm]
\qquad\ds 
\leq~C_2\,h \left(\max_{1\leq k\leq N(j)} \osc\left\{ u(t_j,\cdot)\,;~  \Delta_{jk}\,\right\}
\right)\cdot
\tv\bigl\{ u(t_j,\cdot)\bigr\},
\enda
\eeq
for a suitable constant $C_2$.

Combining (\ref{ee8}) with (\ref{ee9}) and recalling (\ref{tv4}),
we obtain
\bel{es6}
\bega{l}
\ds
\int_{-\infty}^{+\infty} \Big|u(t_{j+1},x) - \bigl( S_{t_{j+1}-t_j} u(t_j,\cdot)\bigr)(x)\Big|\, dx \\[4mm]
\qquad\ds 
\leq~C_3 \left( \kappa_j\,h+  \ve^{2/3}+ h \ve^{1/3}\right)\cdot\sup_{t\in [t_j, t_{j+1}]}
 \tv \bigl\{ u(t,\cdot)\bigr\}.
\enda
\eeq
Recalling that  $h\approx\ve^{1/3}$ and $T=\nu\ve^{1/3}$, from the above analysis we obtain:

\begin{theorem} \label{t:31} Let the basic assumptions {\bf (A1)-(A2)} hold.
Let $t\mapsto u(t,\cdot)\in\D$ be an approximate solution to the Cauchy problem
(\ref{1})-(\ref{2}), taking values in the domain $\D$ of the semigroup
and satisfying {\bf (AL)} and {\bf (P$_\ve$)}.
Then, for some constant $C_4$, the following holds.

Assume that the strip $[0,T]\times \R$ can be covered by 
trapezoids $\Delta_{jk}$, $j=0,\ldots,\nu-1$, $k=1,\ldots N(j)$ as in (\ref{Djk}),
so that (\ref{GG})-(\ref{tv4}) hold.
Then
the difference between $u(T,\cdot)$ and the exact solution 
 $S_T\bar u$ is bounded by
\bel{es9}
  \|u(T,\cdot) - S_T\bar u\|_{\L^1}~\leq~C_4\left(2T+ \sum_{j=0}^{\nu-1}  \kappa_j\right)\ve^{1/3}
 \cdot\sup_{t\in [0,T]}
 \tv \bigl\{ u(t,\cdot)\bigr\}
\,.\eeq
\end{theorem}

{\bf Proof.}
Let $L_0$ be the Lipschitz constant of the semigroup in (\ref{lip2}).
{}From  (\ref{es5}) and (\ref{es6}) it now follows
\bel{es8}\bega{l}
\ds \|u(T,\cdot) - S_T\bar u\|_{\L^1}~\leq~L_0\cdot \sum_{j=0}^{\nu-1}
\Big\| u(t_{j+1}, \cdot) -S_{h} u(t_j, \cdot)\Big\|_{\L^1}\\[4mm]
\qquad \ds \leq ~L_0\cdot \sum_{j=0}^{\nu-1}
C_3 \left( \kappa_j\,h+  \ve^{2/3}+ h \ve^{1/3}\right)\cdot\sup_{t\in [t_j, t_{j+1}]}
 \tv \bigl\{ u(t,\cdot)\bigr\}\\[4mm].
\qquad \ds \leq ~L_0 \, C_3\cdot\left(2T+ \sum_{j=0}^{\nu-1}  \kappa_j\right)\ve^{1/3}
 \cdot\sup_{t\in [0,T]}
 \tv \bigl\{ u(t,\cdot)\bigr\}.
\enda
\eeq
This yields (\ref{es9}), with $C_4 = L_0 C_3$. \endproof

\begin{remark}\label{r:31}{\rm Based on the estimate (\ref{es8}), we seek to understand 
at which rate the error in the approximate solution may approach zero, as $\ve\to 0$.

Having chosen $h\approx \ve^{1/3}$, we can
choose all bounded trapezoids $\Delta_{jk}$, $1< k<N(j)$, to be of diameter $\O(1)\cdot \ve^{1/3}$.
Moreover, by choosing every $b_{j1}$ suitably large and negative, and $a_{j,N(j)}$  large and positive,
we can assume that the solution is nearly constant on the unbounded trapezoids
$\Delta_{j,1}$ and $\Delta_{j, N(j)}$. Here and in  the sequel, the Landau symbol $\O(1)$ denotes a uniformly bounded quantity.

If the exact solution is Lipschitz continuous, we expect 
that the maximum oscillation (\ref{tv4}) will be of size
$\kappa_j = \O(1)\cdot \ve^{1/3}$ for every $j\in\{0,1,\ldots,\nu-1\}$.   
In this case, as $\ve \to 0$ the quantity  $2T+\sum_{j=1}^\nu\kappa_j$ remains uniformly bounded, and the estimate (\ref{es8}) 
indicates that the error vanishes of order $\O(1)\cdot\ve^{1/3}$.

Next, assume that the initial data $\bar u$ contains a jump, generating a centered rarefaction wave
of strength $\sigma$.    In this case, %(see Fig.~\ref{f:hyp222}, left), 
taking into account the decay caused by genuine nonlinearity,
we expect that the oscillation of $u$ over a  trapezoid $\Delta_{jk}$ of diameter $\O(1)\cdot\ve^{1/3}$
will satisfy a bound of the form
\bel{ott}\osc\left\{ u(t,\cdot)\,;~\Delta_{jk}\right\}~= ~\O(1)\cdot\min \left\{ \sigma\,,~{\ve^{1/3}\over  t_j}
\right\}.\eeq
Recalling that $t_j= j\ve^{1/3}$ and $\nu= T\ve^{-1/3}$, this leads to 
\bel{loge}\sum_{j=1}^\nu\kappa_j~=~\sum_{j=1}^\nu \min \left\{ \sigma\,,~{C\ve^{1/3}\over  j\ve^{1/3}}
\right\}~=~\O(1)\cdot \log \nu~=~\O(1)\cdot |\log \ve|.
\eeq
In this case, the estimate (\ref{es8}) would indicate that the error vanishes of order
$\O(1)\cdot \ve^{1/3}\,|\log \ve|$.   The same should hold if the exact solution contains finitely many centered rarefaction waves.

We emphasize, however,  that this is only a heuristic expectation. 
For a numerically computed solution, it needs to be confirmed
by a post-processing algorithm, which can actually provide a bound on the oscillations $\kappa_j$ in 
(\ref{es9}).}
\end{remark}

\section{Solutions with an isolated large shock}
 \label{s:4}
\setcounter{equation}{0}
The error estimates developed in the previous section are not effective for
solutions containing large shocks.  Indeed, around a shock, the oscillation will be large.   As a consequence, even when the diameters of the trapezoids $\Delta_{jk}$ in (\ref{Djk})  approach zero, the maximum oscillation 
$\kappa_j$ in (\ref{tv4}) will remain uniformly large.  For this reason, we do not expect that 
the right hand side of the error bound (\ref{es9}) will approach zero as $\ve\to 0$. 
To cope with this problem,
in this section we develop additional tools to estimate the numerical error in a neighborhood of a shock.

\begin{figure}[ht]
\centerline{\hbox{\includegraphics[width=12cm]{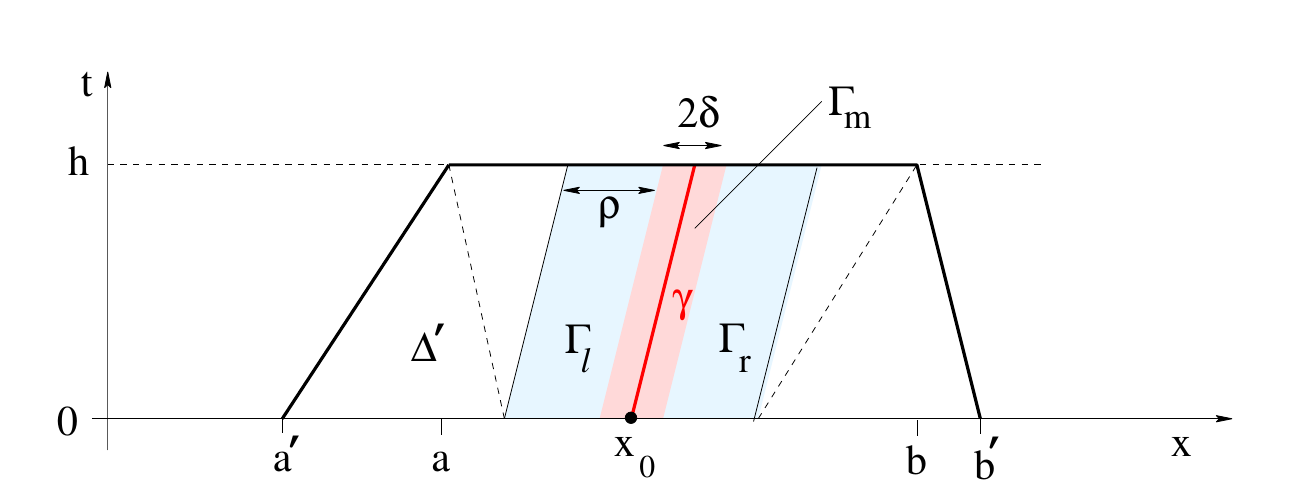}}}
\caption{\small The regions $\Gamma_l,\Gamma_m,\Gamma_r$ introduced at 
(\ref{Glmr}) to trace a large shock, and the trapezoid $\Delta'$  at (\ref{Dp}).}
\label{f:hyp223}
\end{figure}

Consider an approximate solution $u$, which satisfies 
{\bf (AL)} and {\bf (P$_\ve$)}.  We seek a sharper error bound, assuming that the oscillation of 
$u$ is concentrated in a narrow region of the form
\bel{Gdef}
\Gamma~\doteq~\Big\{(t,x)\,;~t\in [0,h],~|x-\gamma(t)|< \delta\Big\},
\qquad\qquad \gamma(t) = x_0 +\lambda t.\eeq
Of course, we expect that $\gamma(\cdot)$ will trace the position of a large shock in the exact solution.
Here $\rho,\delta>0$ are suitable parameters.   Different choices of these values will lead to 
different estimates. As a rule of thumb, it will be useful for the reader 
to keep in mind their order of magnitude:
\bel{oma}
h\,\approx\,\ve^{1/3},\qquad \rho\,\approx\,\ve^{1/3},\qquad \delta\,\approx\, {\ve\over \rho}\,~=~\ve^{2/3}.
\eeq
Referring to Fig.~\ref{f:hyp223}, 
we introduce the points
\bel{aabb}\left\{\bega{rl} a&=~x_0+ \lambda^-h - \delta - \rho,\\[3mm]
b&=~x_0 + \lambda^+h+ \delta + \rho,\enda\right.
\qquad \qquad
\left\{ \bega{rl} a'&\doteq~a-\lambda^+ h,\\[3mm]
 b'&\doteq~b- \lambda^-  h,\enda\right.\eeq
and consider the trapezoidal domain 
\bel{Dp}\bega{rl}
\Delta'&\doteq~\hbox{co}\Big\{(0, a'), \,(0, b'),\, (h,a), \, (h,b)\Big\}\\[3mm]
&=~\Big\{ (t,x)\,;~t\in [0,h]\,,~~a'+\lambda^+t\,\leq\,x\,\leq\, b' +\lambda^-t 
\Big\}.\enda\eeq
Our basic assumption is that, outside the narrow strip $\Gamma$, the oscillation of $u$ remains small.
More precisely, consider the left and right domains
\bel{Dprime}\bega{rl} \Delta'_l&=~
\Big\{ (t,x)\,;~t\in [0,h]\,,~~a'+\lambda^+t\,\leq\,x\,\leq\, x_0 -\delta +\lambda t 
\Big\},\\[3mm]\Delta'_r&=~\Big\{ (t,x)\,;~t\in [0,h]\,,~~x_0 +\delta +\lambda t \,\leq\, x\,\leq\, 
b' +\lambda^- t
\Big\},\enda\eeq
and define
\bel{k'}
\kappa'~\doteq~\max\Big\{\osc\{u\,;~\Delta'_l\}\,,~\osc\{u\,;~\Delta'_r\}\Big\}.\eeq
Calling
\bel{upm}u^-~\doteq~u\bigl(0, x_0-\delta\bigr),\qquad
u^+~\doteq~u\bigl(0, x_0+\delta\bigr),\eeq
the above definition of $\kappa'$ implies
\bel{ues}\bega{rl}
|u(t,x) -u^-| ~\leq ~\kappa'\qquad &\hbox{for}\quad (t,x)\in \Delta'_l\,,\\[3mm]
|u(t,x) -u^+| ~\leq~\kappa'
\qquad &\hbox{for}\quad (t,x)\in \Delta'_r\,.\enda\eeq

Assuming that $\kappa'$ is small, the following result provides a bound on the distance
 between $u$ and the exact solution, 
$\|u(h,\cdot)- S_h\bar u\|_{\L^1([a,b])}$, restricted to the interval $[a,b]$.

\begin{figure}[ht]
\centerline{\hbox{\includegraphics[width=10cm]{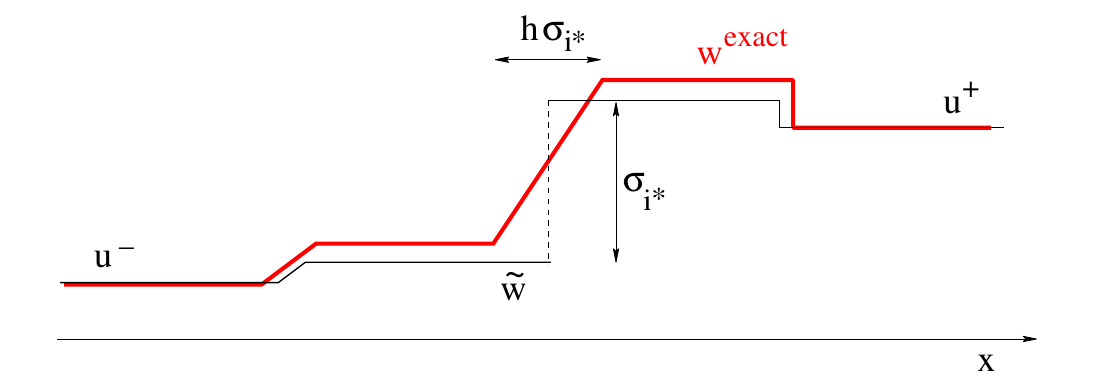}}}
\caption{\small Comparing the entropic solution $w^{exact}$ to the Riemann problem with 
left and right states $u^-,u^+$ with another weak solution $\Tilde w$ containing a non-admissible 
$i^*$-shock of strength $\sigma_{i^*}>0$. Taking into account 
the presence of a centered rarefaction wave in $w^{exact}$, 
The difference between the two solutions can be bounded as
$\|\Tilde w(h,\cdot)- w^{exact}(h,\cdot)\|_{\L^1} =\O(1)\cdot h\sigma_{i^*}^2$.
}
\label{f:hyp226}
\end{figure}

\begin{theorem}\label{t:41} 
Let $t\mapsto u(t,\cdot)\in\D$ be an approximate solution to the Cauchy problem
(\ref{1})-(\ref{2}), taking values in the domain $\D$ of the semigroup,
and satisfying {\bf (AL)} and {\bf (P$_\ve$)}.
Then, for some constant $C_5$, in the above setting we have the error bound
\bel{e41}
\int_a^b\bigl|u(h, x)-(S_h\bar u)(x)\bigr|\, dx~\leq~C_5 \cdot  h \, \left({\ve \over\rho} + \kappa'  + {\rho \kappa' + \delta\over h}\right)^{2/3}
+ C_5\,\Big( \rho\kappa'+ h\kappa' +\delta\Big).\eeq
Moreover, there exists a constant $K_1$ such that, if
\bel{jbig}|u^+-u^-|~\geq~K_1\cdot \left( {\ve\over\delta} + \kappa' + 
{\rho \kappa'+\delta\over h}\right)^{1/3},\eeq
then the estimate (\ref{e41}) can be improved to
\bel{e42}
\int_a^b\bigl|u(h, x)-(S_h\bar u)(x)\bigr|\, dx~\leq~C_5 \cdot  h \, \left({\ve \over\rho} + \kappa'  + {\rho \kappa' + \delta\over h}\right).\eeq

\end{theorem}
\begin{remark}\label{r:41} {\rm It may seem surprising that the error bound (\ref{e42}), valid 
for large jumps, is actually better than (\ref{e41}), which applies to small jumps.   
To understand what is involved here, the following observation can be useful.
If the strength $\sigma = |u^+-u^-|$ is small, it could be that this jump is tracing a
centered rarefaction wave within the exact solution, which gets approximated by a single jump by the
numerical algorithm
(indeed, this is a common feature of front tracking approximations).    
If $\sigma$ is small enough, the entropy produced by the jump is small, and  the assumptions
(\ref{q1})-(\ref{q2}) can still be satisfied.   
This is a ``worst-case scenario":    as shown in Fig.~\ref{f:hyp226},
the corresponding $\L^1$ error is $\O(1)\cdot h \sigma^2$. 
On the other hand, if the strength $\sigma$ of the jump is large, the entropy dissipation
assumption (\ref{q2}) rules out this possibility. 
Therefore, the jump must trace an entropic shock in the exact solution.
}
\end{remark}

{\bf Proof of Theorem \ref{t:41}.} 

{\bf 1.} As  a first step, using {\bf (P$_\ve$)}  we will provide a bound for the error
\bel{E4}
\Big|f(u^+) - f(u^-) - \lambda\bigl( u^+ - u^-\bigr)\Big|.\eeq
As shown in Fig.~\ref{f:hyp223}, denote by $\Gamma_l$, $\Gamma_m$, and $\Gamma_r$ the left, middle, and right domains
\bel{Glmr}\bega{rl}
\Gamma_l&\doteq~\bigl\{(t,x)\,;~t\in  [0,h],~~x\in [\gamma(t)-\delta-\rho\,,~\gamma(t)-\delta]\bigr\},\\[3mm]
\Gamma_m&\doteq~\bigl\{(t,x)\,;~t\in  [0,h],~~x\in [\gamma(t)-\delta\,,~\gamma(t)+\delta]\bigr\},\\[3mm]\Gamma_r&\doteq~\bigl\{(t,x)\,;~t\in  [0,h],~~x\in [\gamma(t)+
\delta\,,~\gamma(t)+\delta+\rho]   \bigr\}.\enda
\eeq
Recalling (\ref{Gdef}), (\ref{aabb}) and (\ref{Dp}), we observe  that the above definitions imply
$$\Gamma_m ~=~\Gamma,\qquad\qquad \Gamma_l\cup\Gamma_m\cup\Gamma_r~\subset ~\Delta'.$$

Given $\rho>0$,
% and a unit vector $\bfw \in\R^n$,
consider a Lipschitz test function $\vp$ such that, for $t\in [0,h]$, one has
\bel{test2}\phi(t,x)~=~\left\{ \bega{cl} 0\qquad &\hbox{if} \quad |x-\gamma(t)|~\geq~\delta+\rho,
\\[3mm]
1 \qquad &\hbox{if} \quad |x-\gamma(t)|~\leq~\delta,
\\[3mm]\ds
{\delta+\rho-|x-\gamma(t)|\over\rho} \qquad &\hbox{if} \quad \delta<|x-\gamma(t)|<\delta+\rho.
\enda\right.\eeq
Then choose any unit vector $\bfw\in\R^n$ and set $\vp(t,x)=\phi(t,x)\, \bfw$.
By construction, for $t\in [0,h]$  the test function 
$\vp$ vanishes outside the union $\Gamma_l\cup\Gamma_m\cup\Gamma_r$.
Notice that 
$$\|\vp_x\|_{\L^\infty} ~= ~{1\over \rho} \,,\qquad\qquad \|\vp_t\|_{\L^\infty} ~= ~{|\lambda|
\over \rho}\,.$$ 
Assuming that the approximate solution $u$ satisfies {\bf ($\bf P_\ve$)}, 
by (\ref{q1}) it follows
\bel{pep}
\bega{l}
\ds\left|\int \vp(0,x) u(0,x)\, dx - \int\vp(h,x) u(h,x)\, dx +
\dint_{ \Gamma_l\cup\Gamma_m\cup\Gamma_r}  u\vp_t+f(u)\vp_x~dxdt 
\right|\\[4mm]
\ds\qquad \leq~C\,\ve\, h  \cdot {\max\bigl\{ 1,\,|\lambda|\bigr\}\over \rho} \cdot \sup_{t\in [0,h]} 
\tv\bigl\{ u(t,\cdot)\bigr\} \,.\enda
\eeq
%{\color{blue}
%$$= ~\O(1)\cdot \ve^{1/2}\cdot h.$$}%upper bound for the oscillation of $u$ on each of the regions
%$$\Delta^-\doteq\,
%\Big\{ (t,x)\,;~t\in [0,h], ~~x\leq x_0-\delta +\lambda^- t\Big\},
%\qquad  \Delta^+\doteq\,\Big\{ (t,x)\,;~t\in [0,h], ~~x\geq
% x_0+\delta +\lambda^+ t\Big\}.$$ 
Using (\ref{ues}), we now estimate
\bel{swk}\left(
\dint_{ \Gamma_l}+\dint_{\Gamma_m}+\dint_{\Gamma_r}  
\right)\bigl\{u\vp_t+f(u)\vp_x\bigr\}~dx\,dt ~\doteq ~
\bfI_l+\bfI_m+\bfI_r \,.
\eeq
Trivially, $\bfI_m=0$ because $\vp_t=\vp_x=0$ on $\Gamma_m$.
By (\ref{ues}) it follows
\bel{Slr} 
\bega{rl}
\ds	\bfI_r+\bfI_l&=\ds~ \dint_{\Gamma_r\cup\Gamma_l}\bigl\{u\vp_t+f(u)\vp_x\bigr\}~dxdt 
\\[4mm]
 %&=~\ds \int_0^h\Big\{\bigl[\lambda u^+-f(u^+) \bigr] -  \bigl[\lambda 
 %u^--f(u^-) \bigr]+\O(1)\cdot\rho\Big\} \bfw \, dt \\[4mm] 
 &=~\ds \int_0^h\Big\{\bigl[\lambda u^+-f(u^+) \bigr] -  \bigl[\lambda 
 u^--f(u^-) \bigr]\Big\} \bfw \, dt +
 \O(1)\cdot  h\,\kappa'.
 \enda
\eeq

Next, by (\ref{k'}) one obtains
\bel{bes}
\bega{l}\ds \left|\int \vp(0,x) u(0,x)\, dx - \int\vp(h,x) u(h,x)\, dx\right|\\[4mm]
\ds
\leq~ \left(\int_{x_0-\delta-\rho}^{x_0-\delta}+ \int_{x_0+\delta}^{x_0+\delta+\rho}\right)
\Big|u(0,x)- u(h, x+\lambda h)\Big|\, dx 
+ \int_{x_0-\delta}^{x_0+\delta} \Big(|u(0,x)| + |u(h,x+\lambda h)|\Big)\, dx\\[4mm]
\ds\leq~\rho\cdot \osc\{ u\,;~\Delta'_l\} + \rho\cdot \osc\{ u\,;~\Delta'_r\} + 
2 \delta\cdot 2\|u\|_{\L^\infty}~\leq~2\rho \kappa' + 4\delta\, \|u\|_{\L^\infty}\,.
%\int_{x_0-\delta-\rho}^{x_0-\delta}\Big(\bigl|u(0,x)-u^-\bigr| 
%+ \bigl|u(h, x+\lambda h)- u^-\bigr|\Big)\,dx\\[4mm]
%\qquad \ds + \int_{x_0+\delta}^{x_0+\delta+\rho}\Big(\bigl|u(0,x)-u^+\bigr| 
%+ \bigl|u(h, x+\lambda h)- u^+\bigr|\Big)\,dx\\[4mm]
%\qquad \ds + \int_{x_0-\delta}^{x_0+\delta} \Big(|u(0,x)| + |u(h,x+\lambda h)|\Big)\, dx\\[4mm]
%=~\O(1)\cdot \rho(K\rho+\kappa)  +\O(1)\cdot \delta\, \|u\|_{\L^\infty}\,.
\enda
\eeq
{}From (\ref{pep}), by (\ref{Slr}) and (\ref{bes}) it follows
\bel{ie1}\bega{l}\ds
\left| \int_0^h\Big\{\lambda ( u^+- u^-) -
\bigl[f(u^+)-f(u^-) \bigr]\Big\} \bfw  \, dt\right|%\\[4mm]\ds\qquad
~=~\O(1)\cdot \left({\ve \, h\over\rho} + h\kappa'  + \rho \kappa' + \delta\right),\enda
\eeq
where the factor $\O(1)$ already accounts for the uniform bound on the total variation.
Choosing the unit vector
$$\bfw ~=~{\lambda ( u^+- u^-) -
\bigl[f(u^+)-f(u^-) \bigr]\over \Big|\lambda ( u^+- u^-) -
\bigl[f(u^+)-f(u^-) \bigr]\Big|}\,,$$
by (\ref{ie1}) we conclude that the error in the Rankine-Hugoniot equations has size
\bel{RHE}\Big|\lambda ( u^+- u^-) -
\bigl[f(u^+)-f(u^-) \bigr]\Big|~=~\O(1)\cdot \left({\ve  h\over\rho} + h\kappa'  +{ \rho \kappa' + \delta
\over h}\right).\eeq
\v
{\bf 2.} Next, consider the piecewise constant function
\bel{pcw}w(t,x)~\doteq~\left\{ \bega{rl} u^-\quad\hbox{if}\quad x< x_0+\lambda t,\\[3mm]
u^+\quad\hbox{if}\quad x> x_0+\lambda t.\enda\right.\eeq
Aim of the next two steps is to prove that 
the difference between $w$ and an exact solution having 
the same initial data is bounded by
\bel{ew}
\|w(h,\cdot)- S_h w(0,\cdot )\|_{\L^1(\R)}~=~\O(1)\cdot h\,\Big|\lambda ( u^+- u^-) -
\bigl[f(u^+)-f(u^-) \bigr]\Big|.\eeq
With this goal in mind, define the averaged Jacobian matrix
$$A~=~\int_0^1 Df\bigl( s u^++(1-s) u^-\bigr)\, ds.$$
Call $\lambda_1<\cdots<\lambda_n$ the eigenvalues of $A$.
Let $\{r_1,\ldots,r_n\}$  and $\{\ell_1,\ldots, \ell_n\}$  be dual bases of right and left eigenvectors of $A$, normalized as (\ref{rli}).   Moreover, let $c_i$, $\ell_{max}$ be such that
\bel{udpm}u^+-u^-~=~\sum_{i=1}^n c_i r_i\,,\qquad\qquad \ell_{max}~\doteq~\max\{|\ell_1|, \ldots,|\ell_n|\}.
\eeq
For every $i=1,\ldots,n$, we then have
\bel{be1}\Big|\ell_i\cdot (\lambda I-A) (u^+-u^-) \Big|~=~|c_i|\, |\lambda-\lambda_i|~\leq~\ell_{max} 
\Big|\lambda ( u^+- u^-) -
\bigl[f(u^+)-f(u^-) \bigr]\Big|.\eeq
  Let $i^*\in \{1,\ldots,n\}$ be a characteristic family such that $|\lambda-\lambda_{i^*}|= \min_i |\lambda-\lambda_i|$. Since the eigenvalues of $A$ are strictly separated, by (\ref{be1}) it follows
\bel{les}|c_i|~=~\O(1)\cdot \Big|\lambda ( u^+- u^-) -
\bigl[f(u^+)-f(u^-) \bigr]\Big|\qquad\forall i\not= i^*.\eeq

We now consider the solution to the Riemann problem with left and right states $u^-, u^+$.
Let $\sigma_1,\ldots,\sigma_n$ be the sizes of the waves in this solution.
As usual, if the $i$-th field is genuinely nonlinear, we choose the sign so that $\sigma_i>0$ 
corresponds to a rarefaction wave, while $\sigma_i<0$ yields an entropy admissible shock.
For future use, we denote by 
\bel{intst}u^-~=~u_0\,,~u_1\,,~\ldots~,~u_n~=~u^+,\eeq
the intermediate states. If the $i^*$-th characteristic field is linearly degenerate, standard estimates on the strength
of these waves yield the bound
\bel{bw1} 
\sum_{i=1}^n |\sigma_i - c_i|~=~\O(1)\cdot \sum_{i\not= i^*} |c_i|\,.\eeq
Indeed, (\ref{bw1}) is trivially true when the right hand side is zero.   The general case is obtained
by an application of the implicit function theorem.  
The same estimate (\ref{bw1}) is achieved when the $i^*$-th field is genuinely nonlinear and $c_{i^*}<0$.
By (\ref{bw1}) it follows
\bel{bw2} |u_{i^*}- u^+| +|u_{i^*-1}-u^-| +\sum_{i\not= i^*} |u_i-u_{i-1}|~=~\O(1)\cdot \sum_{i\not= i^*} |c_i|.
\eeq
In both of the above cases, combining (\ref{be1}), (\ref{les}), and (\ref{bw2}), 
the distance between $w(h,\cdot)$ and an exact solution can be estimated as
\bel{ew1}\bega{l}
\ds{1\over h} \|w(h,\cdot)- S_h w(0,\cdot )\|_{\L^1(\R)}\\[4mm]
\qquad \ds =~\O(1)\cdot \sum_{i\not= i^*} |c_i|
+ \O(1) \cdot |\lambda_{i^*}-\lambda| \,|c_{i^*}|+\O(1)\cdot \bigl(|u_{i^*}- u^+| +|u_{i^*-1}-u^-|\bigr) \, |\lambda|
\\[4mm]
\qquad \ds
=~\O(1)\cdot \Big|\lambda ( u^+- u^-) -
\bigl[f(u^+)-f(u^-) \bigr]\Big|\\[4mm]
\qquad \ds =~\O(1)\cdot \left({\ve \over\rho} + \kappa'  + {\rho \kappa' + \delta\over h}\right).
\enda
\eeq
Notice that the last estimate was obtained from (\ref{RHE}).

\v
{\bf 3.} 
It remains to study the case where the $i^*$-th field is genuinely nonlinear, but $c_{i^*}>0$.   
For this purpose, call $\Tilde w=\Tilde w(t,x)$ the solution to the Riemann problem 
with initial data $\Tilde w(0,\cdot )=w(0, \cdot)$, which contains a non-entropic 
$i^*$-shock of size $\sigma^*$, while all other
waves are entropy admissible.  We observe that all  the above estimates still apply to $\Tilde w$.
In particular, 
\bel{twe}
\ds{1\over h} \|\Tilde w(h,\cdot )-S_h w(0,\cdot)\|_{\L^1(\R)}=~\O(1)\cdot \sum_{i\not= i^*} |c_i|
+ \O(1) \cdot |\lambda_{i^*}-\lambda| \,|c_{i^*}|.\eeq
It remains to estimate the difference between $\Tilde w$ and the entropic solution to the same 
Riemann problem. 
Call $u^-= \Tilde u_0, \Tilde u_1,\ldots, \Tilde u_n= u^+$ 
the intermediate states for the non-entropic solution $\Tilde w$.
 Since shock and rarefaction curves have a second order tangency, comparing with the intermediate states
 (\ref{intst}) of the entropic solution, we find
 \bel{cpi}
 |\Tilde u_i - u_i|~=~\O(1) \cdot \sigma_{i^*}^3\qquad\qquad i=0,1,\ldots,n.\eeq
 Taking into account that the wave connecting the states $ u_{i^*-1}$ 
and $ u_{i^*}$ is a centered rarefaction instead of a single jump, we 
obtain the bound
\bel{wtw}
\ds{1\over h} \|w(h,\cdot)-\Tilde w(h,\cdot )\|_{\L^1(\R)}=~\O(1)\cdot  \sigma_{i^*}^2\,.\eeq
Combining (\ref{twe}) with (\ref{wtw}) we conclude
\bel{tw3}
\ds{1\over h} \|w(h,\cdot )-S_h w(0,\cdot)\|_{\L^1(\R)}=~\O(1)\cdot \sum_{i\not= i^*} |c_i|
+ \O(1) \cdot |\lambda_{i^*}-\lambda| \,|c_{i^*}| + \O(1) \cdot c_{i^*}^2\,.\eeq
%It remains to bound the additional term on the right hand side of (\ref{tw3}).
We claim that the jump $c_{i^*}>0$ must be small, 
otherwise the approximate entropy inequality
(\ref{q2}) would fail. Intuitively, this means that the approximate solution cannot contain a large, 
non-admissible shock.
Indeed,
let $\eta$ be a convex entropy, with entropy flux $q$, such that (\ref{enti}) holds for 
every non-admissible shock of the $i^*$ family.
Let $\phi$ be the test function in (\ref{test2}). Arguing as in (\ref{pep}), by (\ref{q2}) we now 
obtain
\bel{pe3}
\bega{l}
\ds \int \phi(0,x) \eta(u(0,x))\, dx - \int\phi(h,x) \eta(u(h,x))\, dx +\dint_{ \Gamma_l\cup\Gamma_m\cup\Gamma_r}  \eta(u)\phi_t+q(u)\phi_x~dxdt 
\\[4mm]
\ds\qquad \geq~- \O(1)\cdot \,{\ve\over \rho} \, h   \,.\enda
\eeq
As in (\ref{bes}), we have
 \bel{be5} \left|\int \phi(0,x) \eta(u(0,x))\, dx - \int\phi(h,x) \eta(u(h,x))\, dx \right|
~=~\O(1) \cdot (\rho \kappa' + \delta)\,.
\eeq
Repeating the argument at (\ref{swk})-(\ref{Slr}) we obtain
\bel{h1}
\int_0^h\int\bigl\{\eta(u)\phi_t+q(u)\phi_x\bigr\}~dxdt ~=~ h\, \Big\{\bigl[\lambda \eta(u^+)-q(u^+) \bigr] -  \bigl[\lambda 
 \eta(u^-)-q(u^-) \bigr]\Big\} +
 \O(1)\cdot  h\,\kappa'.
\eeq

Next, consider the state $\Tilde u^+$, connected to $u^-$ by a (not entropy admissible)  
$i^*$-shock of size $c_{i^*}>0$.    By the implicit function theorem, one has the bound
\bel{tup}
|\Tilde u^+-u^+|~=~\O(1)\cdot \sum_{i\not= i^*} |c_i|.\eeq
Recalling (\ref{enti}), from (\ref{h1}) we obtain
\bel{h2} \bega{l}
\ds {1\over h} \int_0^h\int\bigl\{\eta(u)\phi_t+q(u)\phi_x\bigr\}~dxdt \\[4mm]
\qquad \ds=~\Big\{\bigl[\lambda \eta(\Tilde u^+)-q(\Tilde u^+) \bigr] -  \bigl[\lambda 
 \eta(u^-)-q(u^-) \bigr]\Big\} +\O(1)\cdot |\Tilde u^+ - u^+| +
 \O(1)\cdot \kappa'\\[4mm]
 \qquad\ds \leq~- c_0|\Tilde u^+- u^-|^3+\O(1)\cdot |\Tilde u^+ - u^+| +
 \O(1)\cdot \kappa' \\[4mm]
 \qquad\ds \leq~-c_0 c_{i^*}^3 + \O(1)\cdot  \sum_{i\not= i^*} |c_i|\\[4mm]
 \ds\qquad \leq~-c_0 c_{i^*}^3 + 
\O(1)\cdot \Big|\lambda ( u^+- u^-) -
\bigl[f(u^+)-f(u^-) \bigr]\Big|, \enda
 \eeq
where (\ref{les}) was used in the last inequality.
Combining (\ref{h2}) with (\ref{pe3}), (\ref{be5}), and (\ref{h1}), and using (\ref{RHE})
to bound the last term in (\ref{h2}), we obtain
\bel{ci3}c_{i^*}~=~
\O(1)\cdot \left({\ve\over\rho} + \kappa'  +{ \rho\over h} \kappa' + {\delta\over h}\right)^{1/3} \,.\eeq
Starting from (\ref{tw3}) and using (\ref{be1}), (\ref{les}), (\ref{RHE}), and (\ref{ci3}), we obtain
\bel{tw4}
\ds{1\over h} \|w(h,\cdot )-S_h w(0,\cdot)\|_{\L^1(\R)}=~\O(1)\cdot \left( {\ve\over\delta} + \kappa' + 
{\rho \kappa'+\delta\over h}\right)^{2/3}.\eeq
\v
{\bf 4.} 
Notice that the estimate (\ref{tw4}) is somewhat weaker, compared with (\ref{ew1}).   In this  step
we show that, if the jump $|u^+-u^-|$ is sufficiently large, then in the genuinely nonlinear case
we must have $c_{i^*}<0$, hence the stronger estimate (\ref{ew1}) holds.
Recalling (\ref{udpm}), notice that
\bel{c**}|c_{i^*}| ~\geq~|u^+-u^-| - \sum_{i\not= i^*} |c_i|~\geq~|u^+-u^-|- \O(1)\cdot \left( {\ve\over\delta} + \kappa' + 
{\rho \kappa'+\delta\over h}\right).\eeq
Therefore, there exists a constant $K_1$ large enough so that, if (\ref{jbig}) holds,
then (\ref{c**}) provides a contradiction with (\ref{ci3}). 
Since (\ref{ci3}) was obtained by assuming that $c_{i^*}>0$, we conclude that 
(\ref{jbig}) is a sufficient condition to guarantee that $c_{i^*}\leq 0$.   In this case, the stronger estimate
(\ref{ew1}) holds.

\v
{\bf 5.} 
Restricted to the interval $[a,b]$, by (\ref{ues}) and (\ref{tw4}),
 the difference between $u$ and the exact solution having 
initial data $u(0,x)=\bar u(x)$ can now be estimated by
\bel{eu1}\bega{l}\ds
\int_a^b\bigl|u(h, x)-(S_h\bar u)(x)\bigr|\, dx\\[4mm]
\qquad\ds
\leq~\left(\int_a^{x_0+\lambda h-\delta} + 
\int_{x_0+\lambda h-\delta}^{x_0+\lambda h+\delta}
+ \int_{x_0+\lambda h+\delta}^b \right)\big| u(h,x)-  w(h,x)\bigr|\, dx\\[4mm]
\qquad\qquad \ds + \|w(h,\cdot)- S_h w(0,\cdot)\|_{\L^1}
+L_0\, \left(\int_{a'}^{x_0-\delta} + \int_{x_0-\delta}^{x_0+\delta}+ \int_{x_0+\delta}^{b'} \right)
 \big| w(0,x)- \bar u(x)\bigr|\, dx
\\[4mm]
\ds
\qquad \leq~\Big[(b-a-2\delta)\,\kappa' + 4\delta \|u\|_{\L^\infty} \Big] 
+ \O(1)\cdot h \, \left({\ve \over\rho} + \kappa'  + {\rho \kappa' + \delta\over h}\right)^{2/3}\\[4mm]
\qquad\qquad +
\Big[(b'-a'-2\delta)\,\kappa' + 4\delta \|u\|_{\L^\infty} \Big] 
 \\[4mm]
\qquad \ds \leq~C_5 \cdot  h \, \left({\ve \over\rho} + \kappa'  + {\rho \kappa' + \delta\over h}\right)^{2/3}
+ C_5( \rho+ h)\kappa' + C_5\delta\,,
\enda
\eeq
for a suitable constant $C_5$.
Indeed, from (\ref{aabb}) it follows
\bel{ba'}b-a~=~2\rho+2\delta + (\lambda^+-\lambda^-)h,\qquad\qquad 
b'-a'~=~2\rho+2\delta + 2(\lambda^+-\lambda^-)h.\eeq
On the other hand, if (\ref{jbig}) holds, then we can use (\ref{ew1}) instead of (\ref{tw4}).
The same argument used in (\ref{eu1})  now yields
\bel{eu2}
\int_a^b\bigl|u(h, x)-(S_h\bar u)(x)\bigr|\, dx~\leq~C_5 \cdot  h \, \left({\ve \over\rho} + \kappa'  + {\rho \kappa' + \delta\over h}\right)+
C_5( \rho+\delta + h)\kappa' + C_5\delta\,.
\eeq

%\begin{remark}\label{r:41} 
%{\rm If $\sigma$ is the size of the (genuinely nonlinear) shock,
%we expect that  for each time $t$, the width of the strip of flagged points
%will be $\delta\approx \O(1) \cdot \ve/\sigma$.   However, assuming
%that  the shock speed is not constant but
%can change in time (at a bounded rate), if we want to enclose the set of 
%flagged points inside a region $\Gamma$ of the form (\ref{Gdef}),
% taking this into account we expect
%$$\delta~=~\O(1)\cdot \left({\ve\over \sigma}+ h^2\right)\,.$$
%If we now choose
%$$\rho~=~\O(1)\cdot \ve^{1/3},\qquad\quad h~=~\O(1)\cdot \ve^{1/3},$$
%then the strength of the shock will be bounded below, say 
%$\sigma\geq c_0\, \rho$.
%Taking this into account, one obtains
%$$\delta~=~\O(1)\cdot \left({\ve\over\sigma}+h^2\right)~=~\O(1)\cdot \ve^{2/3}\,,\qquad {T\over h}~=~\O(1)\cdot \ve^{-1/3}.$$
%Assuming $\kappa=\O(1)\cdot\rho$, and summing the errors in (\ref{eu1}) over $T/h$ steps, we would get a total error
%of size
%$$E~=~\O(1) \cdot {1\over h}\cdot (\ve^{2/3}+\rho h + \rho^2 + \delta)
%~=~\O(1)\cdot \ve^{1/3}. $$
%This is the error generated by a single shock.   Notice that,  if the solution contains finitely many shocks, 
% this leads to an error of the same order of magnitude as in
%(\ref{es14}) which was obtained in the case without large shocks.
%}
%\end{remark}

\section{A post-processing algorithm}
 \label{s:5}
\setcounter{equation}{0}
There are various ways to use the estimates developed in Sections~\ref{s:3} and \ref{s:4}, to obtain a posteriori error bounds.   The underlying idea is to isolate a finite number of thin regions
enclosing the large jumps, where the estimates (\ref{e41}) or (\ref{e42}) can be used.
Then use the bounds (\ref{es9}) on the remaining portion of the domain.

The algorithm described below can be applied to any BV solution of (\ref{1}), but it is designed 
in order to be most effective when the exact solution is  piecewise Lipschitz  with finitely many 
shocks (or contact discontinuities) and centered rarefaction waves.

Let $u:[0,T]\times \R\mapsto \R^n$
be an approximate
solution  of (\ref{1})-(\ref{2}), which satisfies the properties 
{\bf (AL)} and {\bf (P$_\ve$)}.
In this section we introduce an algorithm which checks its total variation,
 identifies the location of large shocks, and constructs trapezoidal subdomains where the 
oscillation remains small.    In view of our previous analysis, this will yield
an  error bound on the $\L^1$ distance (\ref{diff}) between $u$ and an exact
solution.

The algorithm includes three steps.

\v
STEP 1. For each $t\in [0,T]$, we compute
the total variation of $u(t,\cdot)$.   Let $\delta_0>0$ be the constant in 
(\ref{TVd0}).  If
\bel{stv}\sup_{t\in [0,T]}\tv\bigl\{ u(t,\cdot)\bigr\}~\leq~\delta_0,\eeq
then the algorithm can proceed.      On the other hand, if (\ref{stv}) fails, 
the approximate solution may lie outside the domain of the semigroup and no 
error estimate can be provided.    In this case, the algorithm stops.    
\v
STEP 2.  
We now split the interval $[0,T]$ into equal subintervals of size $h= \ve^{1/3}$, 
inserting the times $t_j = j\, h$,  $j=0,1,\ldots,\nu = T/h$.
The next goal is to identify the location of the large shocks,
on each strip $[t_j, t_{j+1}]\times \R$.
For this purpose, we  set $\rho = h = \ve^{1/3}$, $\delta = \ve^{2/3}$, and  
choose two additional parameters:
\begi
\item A lower bound $\sigma_{min}$  for the  size of the jump to be traced.
\item
An upper bound $\kappa'$ for the oscillation of $u$ on a region
to the right and to the left of the jump.
\endi
In view of (\ref{jbig}), it will be convenient to choose these values so that
\bel{jbig2}\sigma_{min}~\geq~K_1\cdot \left( 2\ve^{1/3} + 2\kappa'\right)^{1/3}.\eeq
In this way, the sharper estimate (\ref{e42}) in Theorem~\ref{t:41} will be available.

\begin{figure}[ht]
\centerline{\hbox{\includegraphics[width=12cm]{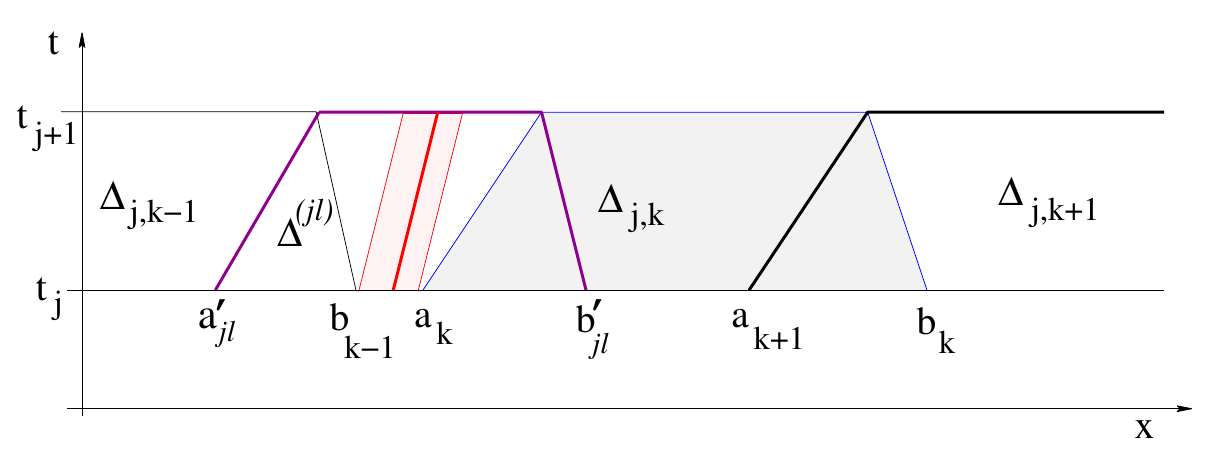}}}
\caption{\small Implementing a post-processing algorithm, each 
strip $[t_j, t_{j+1}]\times\R$ is covered with trapezoids $\Delta_{jk}$ where the 
oscillation remains small (as far as possible), 
and trapezoids $\Delta^{(j\ell)}$ containing a large traced shock.
}
\label{f:hyp224}
\end{figure}

Recalling the construction at (\ref{Gdef})--(\ref{Dprime}), we  introduce
\begin{definition}\label{d:strace}
Given an interval $[x_0-\delta, x_0+\delta]$ and a speed $\lambda\in [\lambda^-,\lambda^+]$, 
consider the polygonal regions
\bel{GJI}\bega{rl} \Gamma&\doteq~\Big\{(t,x)\,;~~ t\in [t_j, t_{j+1}]\,,~~
x_0-\delta+\lambda(t-t_j)\, \leq \,x\,\leq \,x_0-\delta+\lambda(t-t_j)\Big\},\\[3mm]
\Delta'_l&=~
\Big\{ (t,x)\,;~t\in [t_j, t_{j+1}]\,,~~a'+\lambda^+(t-t_j)\,\leq\,x\,\leq\, x_0 -\delta +\lambda (t-t_j) 
\Big\},\\[3mm]\Delta'_r&=~\Big\{ (t,x)\,;~t\in [t_j, t_{j+1}]\,,~~x_0 +\delta +\lambda (t-t_j) \,\leq\, x\,\leq\, 
b' +\lambda^- (t-t_j)
\Big\},\enda\eeq
with $a', b'$ as in (\ref{aabb}).
We say that $\Gamma$ {\bf traces a shock}  during the time interval $[t_i, t_{i+1}]$
if
\bel{osmall}
\max\Big\{\osc\{u\,;~\Delta'_l\}\,,~\osc\{u\,;~\Delta'_r\}\Big\}~\leq~\kappa',\eeq
\bel{jbg}
\bigl| u(t_j,x_0+\delta)- u(t_j, x_0+\delta)\bigr|~\geq~\sigma_{min}\,.\eeq
\end{definition}
In the following, we shall denote by
 \bel{Djl} \Delta^{(j\ell)}~=~\Big\{ (t,x)\,;~t\in [t_j, t_{j+1}]\,,~~a'_{j\ell}+\lambda^+(t-t_j)\,\leq\,x\,\leq\, b'_{j\ell} +\lambda^-( t-t_j) 
\Big\},\qquad \ell=1,\ldots, N'(j),\eeq
the trapezoids containing the traced shocks (see Fig.~\ref{f:hyp224}).

\v
STEP 3.  We cover the remaining region
$[t_j, t_{j+1}]\setminus \bigcup_{\ell=1}^{N'(j)} \Delta^{(j\ell)}$
with finitely many trapezoids of the same form as in (\ref{Djk})
\bel{Tjk}\bega{l}
\Delta_{jk}~=~\bigg\{ (t,x)\,;~t\in [t_j, t_{j+1}],\quad \\[4mm]
\ds\quad {t_{j+1}-t\over h} c_{jk} + {t-t_j\over h} 
(c_{jk} + h \lambda^+ + \ve^{2/3}) ~\leq~x~\leq~{t_{j+1}-t\over h} d_{jk} + {t-t_j\over h} 
(d_{jk} + h \lambda^- - \ve^{2/3}) \bigg\},\enda\eeq
 in such a way that each point $(t,x)\in [t_j, t_{j+1}]\times\R$ 
is contained in at most two of these trapezoids (see Fig.~\ref{f:hyp224}).  
More precisely, we can assume that (\ref{djkk}) holds, for all $k,k'\in \{1,\ldots, N(j)\}$.
Within each time interval $[t_j, t_{j+1}]$,
we compute the maximum oscillation
of $u$ over these trapezoids: 
\bel{kjmax} \kappa_{j}~\doteq~\max_{1\leq k\leq N(j)}~\osc\{ u\,;~\Delta_{jk}\}.\eeq

The next result provides an {\it a posteriori} estimate on the $\L^1$ error in the approximate solution.
Here the estimate refers to the outcome of a
post-processing algorithm, depending on the choice of the parameters $\rho, K$ in 
the definition of flagged points at (\ref{flag}).  We remark that any choice of such parameters 
leads to some error bound.  However, the sharpness of the estimate heavily depends on a suitable 
choice of these parameter values.

\begin{theorem}\label{t:53}
Consider a system of conservation laws satisfying the basic 
assumptions {\bf (A1)-(A2)}. Then there exist constants $C', C''$ such that the following holds.

Let $u:[0,T]\times \R\mapsto\R^n$ be an approximate solution to the Cauchy problem (\ref{1})-(\ref{2}), 
satisfying the conditions {\bf (AL)} and {\bf (P$_\ve$)}, together with (\ref{stv}).
Let $\kappa_j$, $\kappa'$   be the oscillation bounds  in (\ref{tv4}) and (\ref{osmall}),
for a covering with trapezoids $\Delta_{jk}$, $\Delta^{(j\ell)}$  produced by a post-processing algorithm.
Then the difference between $u(T,\cdot)$ and the exact solution is bounded by
\bel{EEE}\bega{l}
\ds \bigl\|u(T,\cdot)- S_T\bar u\bigr\|_{\L^1(\R)}
%\\[4mm]\qquad \ds
~\leq~C'\left( T+ \sum_{j=0}^{\nu-1} \kappa_j
\right) \ve^{1/3}+
C'' \left( \ve^{1/3}\kappa' + \ve^{2/3}\right)\cdot \sum_{j=0}^{\nu-1}N'(j)\,.\enda
\eeq\end{theorem}

{\bf Proof.}  {\bf 1.}
Calling $L_0$ be the Lipschitz constant of the semigroup at (\ref{lip2}), we have
\bel{es0}
\ds \|u(T,\cdot) - S_T\bar u\|_{\L^1}~\leq~L_0\cdot \sum_{j=0}^{\nu-1}
\Big\| u(t_{j+1}, \cdot) -S_h u(t_j, \cdot)\Big\|_{\L^1}\,.\eeq
For each  $j$, in order to estimate the difference $u(t_{j+1}, \cdot) -S_h u(t_j, \cdot)$, 
we consider a covering of the strip $[t_j, t_{j+1}] \times\R$ by 
trapezoids  $\Delta^{(j\ell)}$, $\ell = 1,\ldots N'(j)$ as in (\ref{Djl}), and 
$\Delta_{jk}$, $k=1,\ldots, N(j)$, as in (\ref{Tjk}). 
\v
{\bf 2.}
Recalling (\ref{aabb}), we denote by
$\{t_{j+1}\}\times [a_{j\ell}, b_{j\ell}]$ the upper boundaries of the trapezoids $\Delta^{(j\ell)}$.
These are the trapezoids which contain one large traced shock.
Moreover, we call $\{t_{j+1}\}\times [\Hat c_{jk}, \Hat d_{jk}]$  the upper boundaries
of the remaining trapezoids $\Delta_{jk}$.   
According to  (\ref{Tjk}), this means
$$[\Hat c_{jk}\,,\, \Hat d_{jk}]~=~
\Big[c_{jk} + h \lambda^+ + \ve^{2/3}~,~
d_{jk} + h \lambda^- - \ve^{2/3} \Big].$$
\v
{\bf 3.} 
The same argument used at (\ref{es6}) now yields an error bound on the set 
$$V_j~\doteq~\bigcup_{k=1}^{N(j)}~
[\Hat c_{jk}\,,\, \Hat d_{jk}].$$
Indeed, recalling the uniform bound (\ref{stv}) on the total variation, one obtains
\bel{es66}
\int_{V_j}
 \Big|u(t_{j+1},x) - \left( S_h u(t_j,\cdot)\right)(x)\Big|\, dx ~
\leq~C_3 \left( \kappa_j\,h+  \ve^{2/3}+ h \ve^{1/3}\right)\cdot\delta_0\,.
\eeq
On the other hand, for each $\ell\in\{1,\ldots, N'(j)\}$, applying the estimate
(\ref{e42}) on the interval $ [a_{j\ell}, b_{j\ell}]$
we obtain the error bound
\bel{eu0}
\int_{a_{j\ell}}^{b_{j\ell}} 
\bigl|u(h, x)-(S_h\bar u)(x)\bigr|\, dx~\leq~C_5 \cdot  h \, \left({\ve \over\rho} + \kappa'  + {\rho \kappa' + \delta\over h}\right)
+ C_5\,\left( \rho\kappa'+ h \kappa' +\delta\right).\eeq
\v
{\bf 4.} Recalling our choices 
$$\rho ~=~h~=~\ve^{1/3},\qquad  \delta~=~\ve^{2/3},\qquad  \nu h ~=~\nu \ve^{1/3} ~=~T,$$
summing the terms in (\ref{es66}) over all $j\in\{0,\ldots, \nu-1\}$, 
and summing the terms in (\ref{eu0}) over all $j$ and
all $\ell\in\{ 1,\ldots, N'(j)\}$, we obtain (\ref{EEE}).
\endproof

\begin{remark}\label{r:52} {\rm It is interesting to speculate about the rate at which the error bound on the right hand side of (\ref{EEE}) will approach zero as $\ve\to 0$.
We begin by assuming that the exact solution we are trying to compute is piecewise Lipschitz, with a finite number of centered rarefaction waves, and finitely many non-interacting shocks.

As is Remark~\ref{r:31}, the first term on the right hand side of (\ref{EEE}) is expected to approach zero 
as $\ve^{1/3}|\log\ve|$.   Concerning the second term,  we can fix a constant $C_0$ and 
choose $\kappa' = C_0\ve^{1/3}$.
If the exact solution contains  $N'$ shocks,
we expect that, for all $\ve>0$ sufficiently small, each of these shocks will be traced,
satisfying the inequality (\ref{jbig2}).  
The second term on the right hand side of (\ref{EEE}) will thus have the form
$$C''\left( \ve^{1/3} C_0\ve^{1/3} + \ve^{2/3}\right) \cdot \nu N'~=~\O(1)\cdot \ve^{1/3}.$$
In this case, (\ref{EEE}) would yield
\bel{ERRE} \bigl\|u(T,\cdot)- S_T\bar u\bigr\|_{\L^1(\R)}~=~
\O(1)\cdot \ve^{1/3}|\log\ve| .\eeq
More generally, let us now assume that some of the shocks in the solution interact with each other.
Let $\tau\in [0,T]$ be one of the (finitely many) interaction times.
During a time interval $[\tau^-, \tau^+]$  around $\tau$, of size $\tau^+-\tau^-=\O(1)\cdot \rho$,
we shall not be able to trace the interacting shocks.    As a consequence, 
for $[t_j, t_{j+1}]\cap [\tau^-,\tau^+]\not= \emptyset$, the oscillation 
on one of the  trapezoids $\Delta_{jk}$  in (\ref{Tjk}) (the one which contains a non-traced shock) will be large.  This will force $\kappa_j$ to be large.
However, we expect that the total length of all intervals $[t_j, t_{j+1}]$, where some large shock 
cannot be traced, will have size 
$$\O(1)\cdot \rho\cdot \hbox{[total number of shock interactions]}~=~\O(1)\cdot \ve^{1/3}.$$
In conclusion, the presence of finitely many 
shock interactions will contribute an additional error term $\O(1)\cdot \ve^{1/3}$ to the right hand side of 
(\ref{EEE}).  This will not change its overall order of magnitude.   

One could also argue that, if the solution contains a finite number of compression waves, from which new shocks are formed,  these (non-traced) waves would contribute an error term of the same nature as
a centered rarefaction wave. Therefore, a  bound of the order (\ref{ERRE}) would still
be obtained.

Once again, we emphasize that the bounds (\ref{ERRE}) represent 
only a heuristic expectation. 
For a numerically computed solution, they needs to be confirmed
by a post-processing algorithm, computing a bound on the oscillations $\kappa_j$  in (\ref{tv4}).}
\end{remark}

\section{Properties of  approximation schemes}
\label{s:6}
\setcounter{equation}{0}
In this section we analyze various approximation methods, and check that they
verify the assumptions {\bf (AL)} and {\bf (P$_\ve$)}. Our {\it a posteriori} error estimates can 
thus be applied
to all of them.

\begin{figure}[ht]
\centerline{\hbox{\includegraphics[width=15cm]{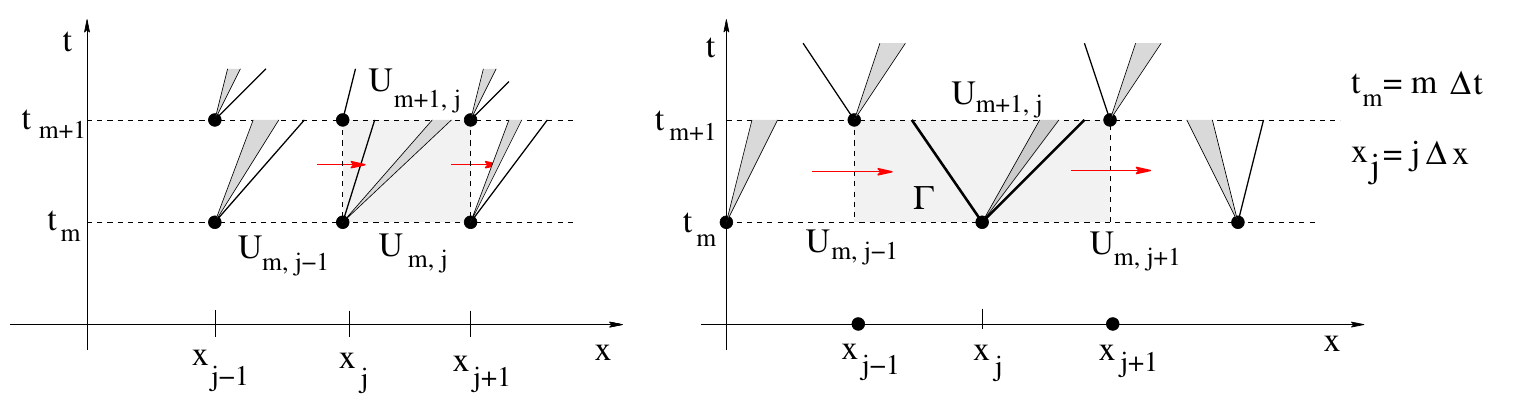}}}
\caption{\small  Left: the upwind Godunov scheme is obtained by solving the Riemann problems at each node $P_{mj}= (t_m, x_j)$,
then by replacing each solution by its average on each of the intervals $[x_j, x_{j+1}]$. 
By the conservation equation, these averages can be 
explicitly computed by (\ref{Ures}). Right: a 
similar construction leads to the Lax-Friedrichs scheme (\ref{LF}).}
\label{f:hyp220}
\end{figure}

\subsection{The Godunov scheme.}

To simplify our discussion, we assume that all characteristic speeds (i.e., all eigenvalues of the
Jacobian matrices $Df(u)$)
lie in the interior of the interval $[0,1]$.  %, for every $u\in\R^n$.
In this case, the Godunov scheme reduces to an upwind scheme.
Given  a mesh size $\ve >0$, consider the grid points 
$$P_{mj} = (t_m, x_j)~=~(\ve m, \ve j).$$
As shown in Fig.~\ref{f:hyp220}, left, 
we consider approximate solutions $u=u(t,x)$ with the following properties:
\begi
\item[(i)] At each time $t=t_m$, the function $u(t_m,\cdot)$ is piecewise constant, namely
$$u(t_m,x)~=~U_{mj}\qquad \hbox{for}\quad x_j< x <x_{j+1}\,,$$
\item[(ii)] For  $t\in [t_m, t_{m+1}[\,$, the function $u(t,\cdot)$ yields the exact 
solution to (\ref{1}) with initial data $u(t_m, \cdot)$.    This is obtained by solving the Riemann problems at each node $x_j$. 
\item[(iii)] At time $t_{m+1}$, we take the average of 
$u(t_{m+1}-\,, \cdot)$ on each interval $[x_j, x_{j+1}]$. Namely
\bel{ave}
u(t_{m+1}, x)~=~U_{m+1, j}~\doteq~
{1\over\ve} \int_{x_j}^{x_{j+1}} u(t_{m+1}, y)\, dy\qquad 
\hbox{for}\quad x_j<x<x_{j+1}\,.\eeq
Since we are assuming that all wave speeds are contained in the interval $[0,1]$, using the conservation equations
these average values $U_{m+1,j}$ can be computed by 
\bel{Ures} U_{m+1,\,j}~=~U_{m,j}+\Big(
f( U_{m,j-1})-f(U_{m,j})\Big).\eeq
\endi
We check that 
an approximate solution $u$ produced by the Godunov scheme with mesh size $\ve>0$ satisfies the Lipschitz condition {\bf (AL)}.  Indeed,
for every $0\leq \tau<\tau' \leq T$ with $\tau,\tau'\in \ve\N$,
\bel{lipg}\bega{l}\bigl\| u(\tau',\cdot)-u(\tau,\cdot)\bigr\|_{\L^1(\R)}~\leq 
\ds\sum_{\tau<t_m\leq\tau'}\sum_{j}\int_{x_{j-1}}^{x_j}\bigl|u(t_{m+1},x)-u(t_m,x)\bigr|\,dx\\[4mm]
\ds	\qquad=~\sum_{\tau<t_m\leq\tau'}\sum_j \ve\Big|
f( U_{m,j-1})-f(U_{m,j})\Big|~\leq~\ve\,\sum_{\tau<t_m\leq\tau'}\hbox{Tot.Var.} \Big\{f(u(t_m,\cdot))\Big\}
\\[4mm]
\ds\qquad\leq~
(\tau'-\tau) \cdot \hbox{Lip}(f)\cdot \sup_{t\in[\tau,\tau']} \tv\bigl\{u(t,\cdot) \bigr\}.
\enda
\eeq
\v

To prove that the property {\bf (P$_\ve$)} also holds, we shall use
\begin{lemma}\label{l:61}
Let $w:[0,\ve]\mapsto\R$ be any function with bounded variation, and assume
$\vp\in \C^1$.  Consider the average value
$$\ov w~\doteq~{1\over\ve} \int_0^\ve w(y)\, dy.$$
Then 
\bel{e6} \left|\int_0^\ve[w(x)-\ov w]\,\vp(x)\, dx\right|~\leq~\tv\bigl\{w\,;~]0,\ve[\,\bigr\}
\cdot \ve^2 \|\vp_x\|_{\L^\infty}\,.\eeq
\end{lemma}

{\bf Proof.} Call $\ov\vp$ the average value of $\vp$ over $[0,\ve]$.  
Then 
\bel{e7} \bega{l} \ds\left|\int_0^\ve[w(x)-\ov w]\,\vp(x)\, dx\right|~
=~\left|\int_0^\ve[w(x)-\ov w]\,
(\vp(x)-\ov \vp)\, dx\right|\\[4mm]
\qquad\leq\ds \int_0^\ve \|w(\cdot)-\ov w\|_{\L^\infty}\cdot \,
\|\vp(\cdot )-\ov \vp\|_{\L^\infty} \, dx
~\leq~\ve\cdot \tv\bigl\{w\,;~]0,\ve[\,\bigr\}\cdot \ve \|\vp_x\|_{\L^\infty}\,.
\enda\eeq\endproof
\v
Next, fix $\ve>0$ and consider any test function $\vp\in \C^1_c(\R^2)$.
Since the Godunov approximations coincide with exact solutions on each of the
half-open intervals $[t_m, t_{m+1}[\,$, we have 
\bel{e5} \bega{l}\ds \left| \int u(\tau,x)\vp(\tau,x)\, dx-\int u(\tau',x)\vp(\tau',x)\, dx+
\int_{\tau}^{\tau'}\int \bigl\{u\vp_t+f(u)\vp_x\bigr\}~ dxdt\right|\\[4mm]
\qquad \ds =~\left|\sum_{\tau <t_m\leq \tau'} \int [u (t_m, x) - u(t_m-, x)]\vp(t_m,x)\, dx\right|.
\enda
\eeq
Using Lemma~\ref{l:61} we obtain
\bel{e8}\bega{l}\ds\sum_{\tau<t_m\leq \tau'}
\left| \int [u (t_m, x) - u(t_m-, x)]\,\vp(t_m,x)\, dx\right|\\[4mm]
\ds\qquad \leq~
\sum_{\tau<t_m\leq \tau'} \sum_j \ve^2 \, \tv\Big\{u(t_m,\cdot) \,;~]j\ve,(j+1)\ve[\,\Big\}
\cdot  \|\vp_x(t_m,\cdot)\|_{\L^\infty}\\[5mm]
\qquad\ds
=~\ve\,(\tau'-\tau) \|\vp_x\|_{\L^\infty}\cdot \sup_{t\in[\tau,\tau']} \tv\bigl\{u(t,\cdot) \bigr\}.
%\leq ~\ve ( \ve+T) \cdot \,\sup_{t\in [0,T]} \tv\bigl\{u(t,\cdot) \bigr\}\cdot  \|\vp_x\|_{\L^\infty}\,.
\enda
\eeq
This yields (\ref{q1}).
In order to prove {\bf (P$_\ve$)}, given a convex entropy with entropy flux $q$,
it remains to check that (\ref{q2}) is satisfied as well.
Let $\varphi\geq 0$ be a test function in $\mathcal{C}^1_c(\R^2)$. Integration by parts yields
\bel{eg}\bega{l}
\ds\int \eta(u(\tau,x))\vp(\tau,x)\, dx-\int \eta(u(\tau',x))\vp(\tau',x)\, dx+	\int_{\tau}^{\tau'}\int{\eta(u)\varphi_t + q(u)\varphi_x}~dt~dx \\[3mm]
\ds	=~-\int_{\tau}^{\tau'}\int \lbrace\eta(u)_t+q(u)_x\rbrace\varphi~dx~dt+\sum_{\tau<t_m\leq \tau'} \int(\eta(u(t_m-,x))-\eta(u(t_m,x)))\varphi(t_m,x)dx. \enda\eeq
By construction, the approximation $u$ is an entropy weak solution of the hyperbolic system of conservation law  in every strip $[t_{m-1},t_m[\times \R$, therefore the first term on the right hand side 
of (\ref{eg}) is non negative.

By the convexity of $\eta$, we can apply Jensen's inequality and obtain
\bel{J}
	\eta(u(t_m,x))~= ~\eta\Big( {1\over \ve}\int_{x_j}^{x_{j+1}} u(t_{m}-, y)\, dy\Big)~\leq~
	 {1\over \ve}\int_{x_j}^{x_{j+1}} \eta(u(t_{m}-, y))\, dy\, ,
\eeq
for $x_j<x<x_{j+1}$. In turn, this yields
\bel{m}\bega{l}
	\ds \sum_{\tau <t_m\leq \tau'} \int_{\R}(\eta(u(t_m-,x))-\eta(u(t_m,x)))\varphi(t_m,x)dx \\[3mm] 
   \ds~	= ~\sum_{\tau <t_m\leq \tau'}\sum_j \int_{x_j}^{x_{j+1}}(\eta(u(t_m-,x))-\eta(u(t_m,x)))\varphi(t_m,x)\, dx \\[3mm] 
	\ds ~ \geq  ~ \sum_{\tau <t_m\leq \tau'}\sum_j \int_{x_j}^{x_{j+1}}\Big(\eta(u(t_m-,x))-\int_{x_j}^{x_{j+1}} \eta(u(t_{m}-, y))\, dy\Big)\varphi(t_m,x)\, dx\\[3mm] 
	\ds ~ \geq  ~ -\, \hbox{Lip}(\eta)\,\ve\cdot (\tau'-\tau) \|\vp_x\|_{\L^\infty}\cdot \sup_{t\in [\tau,\tau']} \tv\bigl\{u(t,\cdot) \bigr\}.
\enda\eeq

\subsection{The Lax-Friedrichs scheme.}
Consider  step sizes $ \Delta t, \Delta x>0$
so that all characteristic speeds satisfy the CFL condition
\bel{CFL} |\lambda_i|~<~{\Delta x\over \Delta t}\,.\eeq 
As shown in Fig.~\ref{f:hyp220}, right, we then construct a staggered grid
with nodes at the points 
$$P_{mj}~=~(m\,\Delta t, \, j\, \Delta x),\qquad\qquad m+j~\hbox{even}.$$
The Lax-Friedrichs approximations are  defined inductively as follows.
Given a piecewise constant function 
$u(t_m,\cdot)$, with jumps at the points $P_{mj}$ with $m+j$ even,
for $t\in [t_m, t_{m+1}[$ we let $u(t,\cdot)$ be the exact solution of the 
system of conservation laws (\ref{1}) with the given data at $t=t_m$.
We then define $u(t_{m+1}, \cdot)$ to be the piecewise constant function
obtained by taking the average of $u(t_{m+1}-,\cdot)$ over every
interval $[x_{j-1}, x_{j+1}]$ with $m+j$ even.
By the conservation equations,  if all characteristic speeds 
satisfy $|\lambda_i|< \Delta x/\Delta t$, these average values 
$$U_{m+1,j}~=~{1\over 2\Delta x}\int_{x_{j-1}}^{x_{j+1}} u(t_{m+1}-,\, x)\, dx
\,,\qquad\qquad
m+j ~\hbox{even},$$ are
inductively computed by
the Lax-Friedrichs scheme
\bel{LF}U_{m+1,j}~=~{1\over 2}(U_{m,j+1}+ U_{m, j-1}) -{\Delta t\over 2\Delta x}\bigl[f(U_{m, j+1}) -
f(U_{m, j-1})\bigr].\eeq
Setting $\ve = \Delta t$, both the approximate Lipschitz condition 
{\bf (AL)} and the property {\bf (P$_\ve$)} can be proved in the same 
way as for the Godunov scheme.  We thus omit details.

\subsection{Backward Euler approximations.}

We now discretize time but keep  space continuous. 
We assume that all characteristic speeds  are strictly positive.
Calling $\ve=\Delta t$ the time step, and 
setting $U_m(x)= u(m\ve,x)$, the backward Euler approximations are defined in terms of 
implicit equations
$$u(t+\ve,x)~=~u(t,x)-\ve f(u(t+\ve,x))_x\,.$$
Equivalently: 
\bel{BEA}
U_{m+1}(x)~=~U_m(x) - \ve \, f(U_{m+1}(x))_x\,.\eeq
At the present  time, the convergence of these approximations for general $n\times n$ hyperbolic systems
is not known.   Complete results are available in the scalar case \cite{Crandall}, which can be
handled by the 
general theory of nonlinear contractive semigroups \cite{CL}.
Let $U_m(\cdot)$ be a sequence of solutions to (\ref{BEA}) with $m=0,1,\ldots$
 and define the approximate solution $u$ by setting 
\[u (t,x)=U_m (x) \qquad \hbox{ for }\quad m\ve\leq t <(m+1)\ve.\]
Then \textbf{(AL)} follows by 
$$ \bega{rl}
\bigl\| u(\tau,\cdot)-u(\tau',\cdot)\bigr\|_{\L^1(\R)}&\leq~\ds\sum_{\tau <t_m\leq \tau'}\bigl\| u(t_m,\cdot)-u(t_{m-1},\cdot)\bigr\|_{\L^1(\R)}\\[4mm]
&=~ \ds\ve \sum_{\tau <t_m\leq \tau'}\int \bigl| f(u(t_m, x))_x\bigr|\, dx\\[4mm]
&\leq ~\hbox{Lip}(f)\,(\tau'-\tau)\cdot \Big(\sup_{t\in[\tau,\tau']} \hbox{Tot.Var.}\lbrace u(t, \cdot)\rbrace \Big).
\enda $$
Next, we check that {\bf (P$_\ve$)} holds.
As before, fix $0\leq \tau<\tau'\leq T$ with $\tau,\tau'\in\ve\N$. Given  a test function
$\varphi\in\mathcal{C}^1_c(\R^2)$,  we compute
$$ \bega{l}\ds
	\left|\int u(\tau,x)\vp(\tau,x)\, dx-\int u(\tau',x)\vp(\tau',x)\, dx+\int_{\tau}^{\tau'}\int\bigl\{u\varphi_t + f(u)\varphi_x\bigr\}\,dx\,dt\right|\\[4mm]
\qquad	=\ds~\left| -\int_{\tau}^{\tau'}\int f(u)_x \varphi \,dx\,dt + \int_\R
	{\sum_{\tau <t_m\leq \tau'} \Big(u(t_{m-1},x)-u(t_m,x)\Big)\varphi(t_{m-1},x)}~dx\right| \\[4mm]
\qquad		\ds \leq~\bigg|-\int \left(\sum_{\tau <t_m\leq \tau'}\frac{u(t_{m-1},x)-u(t_m,x)}{\ve}\int_{t_{m-1}}^{t_m}\varphi(t,x)dt \right) dx\\[4mm]
	\qquad\qquad\qquad \ds	+ \int\sum_{\tau <t_m\leq \tau'}{ \Big(u(t_{m-1},x)-u(t_m,x)\Big)\varphi(t_{m-1},x)}\,dx\bigg| \\[4mm]
\qquad		\ds\leq~\int \sum_{\tau <t_m\leq \tau'} \big| u(t_{m-1},x)-u(t_m,x) \big| \left|\varphi(t_{m-1},x)- \frac{1}{\ve}\int_{t_{m-1}}^{t_m}\varphi(t,x)dt \right|dx\\[4mm] 
%	\ds\leq~ \ve \Vert \varphi_t \Vert_{\L^\infty} \sum_{m=1}^{T/\ve}\int_{\R}\big| u(t_{m-1},x)-u(t_m,x) \big|dx\,.
\qquad		\ds \leq ~ \ve\,(\tau'-\tau)\,\Vert \varphi \Vert_{W^{1,\infty}} \hbox{Lip}(f) \cdot \Big(\sup_{t\in[\tau,\tau']} \hbox{Tot.Var.}\lbrace u(t, \cdot)\rbrace \Big).	
	\enda $$	
	
Finally, let $\eta$ be a convex entropy with entropy flux $q$. If $\varphi \geq 0$ is a test function in $\mathcal{C}^1_c(\R^2)$, we have 
	$$ \bega{l}
	\ds\int \eta(u(\tau,x))\vp(\tau,x)\, dx-\int \eta(u(\tau',x))\vp(\tau',x)\, dx+\int_{\tau}^{\tau'}\int\lbrace\eta(u)\varphi_t + q(u)\varphi_x\rbrace~dx~dt\\[4mm]
	~=~\ds - \int_{\tau}^{\tau'}\int q(u)_x \varphi ~dt~dx + \int{\sum_{\tau<t_m\leq \tau'} \Big(\eta(u(t_{m-1},x))-\eta(u(t_m,x))\Big)					\varphi(t_{m-1},x)}~dx \\[4mm]
	~=~ \ds- \sum_{\tau <t_m\leq \tau'}\int \Big( D\eta(u(t_m,x))\cdot\frac{u(t_{m-1},x)-u(t_m,x)}{\ve}\int_{t_{m-1}}^{t_m}\varphi(t,x)dt \Big) dx\\[4mm]
	\qquad \qquad\qquad\qquad\ds+\int{\sum_{\tau <t_m\leq \tau'}  \Big(\eta(u(t_{m-1},x))-\eta(u(t_m,x))\Big)\varphi(t_{m-1},x)}\,dx \\[4mm]
	\ds\geq ~\sum_{\tau <t_m\leq \tau'} \int\Big(\eta(u(t_{m-1},x))-\eta(u(t_m,x)\Big)\Big(\varphi(t_{m-1},x)- \frac{1}{\ve}\int_{t_{m-1}}^{t_m}\varphi(t,x)dt \Big)dx\\[4mm]
	\ds\geq ~-\ve\,(\tau'-\tau)\,\Vert \varphi \Vert_{W^{1,\infty}} \hbox{Lip}(\eta)\hbox{Lip}(f) \cdot\Big(\sup_{t\in[\tau,\tau']} \hbox{Tot.Var.}\lbrace u(t, \cdot)\rbrace \Big).
	\enda $$
Notice that the convexity of $\eta$ was here used in the inequality 
	\bel{ucx}	D\eta(u(a,x))\cdot \Big(u(a,x)-u(b,x)\Big)~\geq ~\eta(u(a,x))-\eta(u(b,x).
	\eeq
\v
\subsection{The smoothing method.}
Next, we  consider an approximate solution to (\ref{1})-(\ref{2}) obtained by 
periodic mollifications, taking the convolution with a smoothing kernel
$K\in \C^\infty_c(\R)$.  As usual, we assume 
\bel{Kr}\left\{
\bega{rl} K(x)&> ~0\quad\hbox{for}~ |x|<1,\\[2mm]
K(x)&=~0 \quad\hbox{for}~ |x|\geq 1,
\enda\right.\qquad K(x)=K(-x),\qquad \int K(x)\, dx~=~1,\eeq
and set $K_\delta(x)\doteq \delta^{-1} K(\delta^{-1}x)$.
We fix a time step $\ve>0$ and define an $\ve$-approximate solution $u$ by setting
$$t_m~=~m\,\ve, \qquad u(t_m ,\cdot)~=~K_\delta * u (t_m-,\cdot),$$
and letting $u$ be a classical solution to (\ref{1}) on each 
half-open interval $[t_m, t_{m+1}[\,$.

As in the scalar case (see \cite{HR}), the method is well-defined provided that the ratio
$\ve/\delta$ is suitably small.
To see this, in connection with the quasilinear system
\bel{quas1}u_t+A(u) u_x~=~0,\qquad\qquad A(u)\,=\, Df(u),\eeq
we choose bases $\{r_1(u),\ldots, r_n(u)\}$ and 
$\{l_1(u),\ldots, l_n(u)\}$ of right and left eigenvectors of $A(u)$, 
normalized so that 
\bel{dualb} |r_i(u)|\,=\,0,\qquad\qquad l_i(u)\cdot r_i(u)~=~\left\{\bega{rl} 1\quad\hbox{if}~i=j,\cr 0\quad\hbox{if}~i\not= j.\enda\right.\eeq
We denote by 
$u_x^i\,=\, l_i\cdot u_x$
the $i$-th component of the gradient vector $u_x$ w.r.t.~this basis.
{}From (\ref{dualb}) and  (\ref{quas1}) it follows
$$u_x~ =~\sum_{i=1}^n u_x^i r_i(u)\qquad\qquad  u_t ~= ~
- \sum_{i=1}^n \lambda_i(u) u^i_x r_i(u).$$
Differentiating the first equation w.r.t.~$t$ and the second one w.r.t.~$x$,
then equating the results, one obtains
a semilinear 
system of evolution equations for the scalar components $u_x^i$, having 
the form
\bel{evcomp}
(u_x^i)_t + \lambda_i (u_x^i)_x ~=~\sum_{j,k=1}^n
g^i_{jk}(u)u_x^j u_x^k.\eeq
See for example Section 1.6 in \cite{Btut} for details.
Assume that
$$\bigl| g^i_{jk}(u)\bigr|~<~M_g$$
for all $i,j,k$, and all $u$ in the domain were the solution is defined.
Let $t\mapsto Z(t)$ be the solution to the ODE
\bel{odez}
{d\over dt} Z(t)~=~n^2 M_g \, Z^2(t),\qquad\qquad Z(0)\,=\, Z_0\,.\eeq
Assume that, at time $t=0$,  there holds
\bel{iode}\bigl|u^i_x(0,x)\bigr|~\leq~Z_0\qquad\forall x\in\R, ~~i=1,\ldots,n.\eeq
A comparison argument now yields
$$\bigl|u^i_x(t,x)\bigr|~\leq~Z(t)\qquad \forall x\in\R, ~t\in [0, T_0[\,,$$
where 
$$T_0~=~ {1\over n^2 M_g \, Z_0 }$$
is the time where the solution to (\ref{odez}) blows up.\\
It remains to give an upper bound for the gradient components
after each mollification. This is achieved  observing that  
$$\bega{l}\ds\|u_x^i(t_m,\cdot)\|_{\L^\infty}~\leq~ \|u_x(t_m,
\cdot)\|_{\L^\infty}\cdot \sup_u \,|l_i(u)|\\[3mm]
\qquad ~\leq~ \ds{\tv\lbrace u(t_{m}-, \cdot)\rbrace\over \delta}\,\Vert K 
\Vert_{\L^{\infty}}\cdot \sup_u \,|l_i(u)|.\enda$$
Therefore, if we choose
$$ 0<\ve< {\delta \over n^2\,M_g\, \Vert K \Vert_{\L^{\infty}}}\,\Big(\sup_m \tv\lbrace u(t_{m}-, \cdot)\rbrace\Big)^{-1} \Big(\sup_i\sup_u |l_i(u)|\Big)^{-1}, $$
all the components $u^i_x$ remain bounded on each strip $[t_m,t_{m+1}[\times \R$, and the
approximate solution is well-defined.
%For instance we can assume 
%$$ \delta = 2\, n^2\,M_g\,\Vert K \Vert_{\L^{\infty}}\,\Big(\sup_m \tv\lbrace u(t_{m}-, \cdot)\rbrace\Big) \Big(\sup_i\sup_u |l_i(u)|\Big)\, \ve.$$

We now check that the assumption {\bf (AL)} holds:
$$ \bega{lr}
\bigl\| u(\tau',\cdot)-u(\tau,\cdot)\bigr\|_{\L^1(\R)}\\[4mm]
\ds~\leq ~\sum_{\tau<t_m\leq\tau'}\int \bigl|~K_\delta * u (t_{m}-,x)-u(t_{m}-,x)|\, dx + \sum_{\tau<t_m\leq\tau'}\int \bigl|u (t_{m}-,x)-u(t_{m-1},x)|\, dx\\[4mm]
\ds~\leq ~\delta\Big({\tau'-\tau\over \ve}\Big) \cdot\sup_{t\in[\tau,\tau']} \tv\bigl\{ u(t,\cdot)\bigr\} + L\,(\tau'-\tau)\cdot\sup_{t\in[\tau,\tau']} \tv\bigl\{ u(t,\cdot)\bigr\}\\[4mm]
%\ds~\leq~ (1+L)({b-a\over \ve}+1)\sup_{t_j\leq t\leq t_{j+1}} \tv\bigl\{ u(t,\cdot)\bigr\},
\ds~ = ~ C\,(\tau'-\tau)\cdot\sup_{t\in[\tau,\tau']} \tv\bigl\{ u(t,\cdot)\bigr\}.
\enda $$
Here we are using the  estimate
\bel{br} \bega{lr}
\ds\int \bigl|~K_\delta * u (t_{m}-,x)-u(t_{m}-,x)|\, dx\\[4mm]
~ \ds\leq ~\int_{-\delta}^{\delta}K_\delta(y) \sup_{|y|\leq \ve}\Big(\int |u(t_{m}-,x)-u (t_{m}-,x-y)|\, dx\Big)\,dy\\[4mm]
~\leq ~\delta ~\tv\lbrace u(t_{m}-, \cdot)\rbrace.
\enda \eeq

To prove {\bf (P$_\ve$)} we shall need the following result.

\begin{lemma}\label{l:62}
	Let $w:\R\mapsto\R$ be any function with bounded variation, and assume
	$\vp\in \C^1$.  Let  $K\in\C^{\infty}_c$ be a smoothing kernel as in (\ref{Kr}) and define
	$$ \Tilde{w}~\doteq~ K_\delta\ast w.$$
Then 
	\bel{ls} \left|\int [\Tilde{w}(x)-w(x)]\,\vp(x)\, dx\right|~=~\O(1)\cdot \delta^2\cdot \|\vp\|_{W^{1,\infty}}\cdot \tv\{w\}\,.\eeq
\end{lemma}
{\bf Proof.}
We rewrite the left hand side of (\ref{ls}) in a more suitable way:
\bel{fi}\bega{l} 
\ds\left|\int [\Tilde{w}(x)-w(x)]\,\vp(x)\, dx\right|\\[3mm]
~=~\ds\left|\dint K_\delta(x-y) [w(y)-w(x)] \,\vp(x)\, dy\,dx\right|\\[3mm]
~=~\ds\left|\dint K_\delta(x-y) [\vp(y)-\vp(x)] \,w(x)\, dy\,dx\right|\\[3mm]
~=~\ds\left|\int \left(\int_{-\infty}^x \int  K_\delta(x-y) [\vp(y)-\vp(x)] \, dy dz \right) w'(x)\, dx\right|.
\enda \eeq 
Next we prove the estimate
\bel{IK}
\left|\int_{-\infty}^x\int K_\delta(z-y) [\vp(y)-\vp(z)] \, dy dz\right| ~=~\O(1)\cdot \delta^2\cdot \|\vp\|_{W^{1,\infty}}\,.\eeq

The integral in (\ref{IK}) can be split in
\bel{spl} \bega{l}\ds \int_{-\infty}^x\int K_\delta(z-y) [\vp(y)-\vp(z)] \, dy dz\\[3mm] ~=~\ds\int_{-\infty}^{x-\delta}\int_{z-\delta}^{z+\delta} K_\delta(z-y) [\vp(y)-\vp(z)] \, dy dz + \int_{x-\delta}^{x}\int_{z-\delta}^{z+2x-\delta} K_\delta(z-y) [\vp(y)-\vp(z)] \, dy dz \\[3mm]
\quad \quad \quad \quad~+~\ds\int_{x-\delta}^{x}\int_{z+2x-\delta}^{z+\delta} K_\delta(z-y) [\vp(y)-\vp(z)] \, dy dz 
\enda\eeq
and by introducing the change of coordinates $y={q-t\over 2}$ and $z={q+t\over 2}$ it follows that 

$$\bega{l}
\ds\int_{-\infty}^{x-\delta}\int_{z-\delta}^{z+\delta} K_\delta(z-y) [\vp(y)-\vp(z)] \, dy dz + \int_{x-\delta}^{x}\int_{z-\delta}^{z+2x-\delta} K_\delta(z-y) [\vp(y)-\vp(z)] \, dy dz \\[3mm]
\quad \quad~=~\ds \int_{-\infty}^{2x-\delta}\int_{-\delta}^{\delta}K_\delta(t)\Big[\vp\Big({q-t\over 2}\Big)-\vp({q+t\over 2})\Big]dt\,dq\\[3mm]
\quad \quad~=~\ds 0.
\enda $$
Therefore the only contribution is given by the last integral in (\ref{spl}). Calling 
$\Sigma$ be its domain of integration, we find  
\bel{IK3}\bega{rl}
\ds\left|\dint_\Sigma K_\delta(z-y) [\vp(y)-\vp(z)] \, dy dz \right|
&\ds \leq~ \meas(\Sigma)\cdot \|K_\delta\|_{\L^\infty}\cdot \sup_{(x,y)\in \Sigma} 
|\vp(y)-\vp(z)|\\[4mm]
&\ds =~\delta^2 \cdot C\delta^{-1} \cdot 2\delta \, \|\vp\|_{W^{1,\infty}}\,,
\enda
\eeq
Together with (\ref{fi}), this yields (\ref{ls}).
\endproof
\v
Now consider any test function $\varphi\in\mathcal{C}^1_c(\R^2)$.  Using the above lemma, we obtain
$$ \bega{lr}\ds
\left|\int u(\tau,x)\vp(\tau,x)\, dx-\int u(\tau',x)\vp(\tau',x)\, dx+\int_{\tau}^{\tau'}\int\bigl\{u\varphi_t + f(u)\varphi_x\bigr\}\,dx\,dt \right|\\[4mm]
~\ds=~\left|\int\sum_{\tau< t_j\leq\tau'} \bigl(u (t_{j},x)-u(t_{j}-,x)\bigr)\,\varphi(t_{j},x) dx\right|\\[4mm]
~\ds\leq~ C_1\,\delta^2\Big({\tau'-\tau\over \ve}\Big)\Vert \varphi \Vert_{W^{1,\infty}}\cdot \Big(\sup_{t\in[\tau,\tau']} \hbox{Tot.Var.}\lbrace u(t, \cdot)\rbrace \Big)\\[4mm]
~\ds=~ C_2\, \ve \,(\tau'-\tau)\Vert \varphi\Vert_{W^{1,\infty}}\cdot\Big(\sup_{t\in[\tau,\tau']} \hbox{Tot.Var.}\lbrace u(t, \cdot)\rbrace \Big).
\enda$$

Finally, let $\eta$ be a convex entropy with entropy flux $q$. 
For any non-negative test function $\vp\in \C^1_c(\R^2)$  one has
$$ \bega{lr}\ds
\int\eta(u(\tau,x))\vp(\tau,x)\, dx-\int \eta(u(\tau',x))\vp(\tau',x)\, dx+\int_{\tau}^{\tau'}\int\bigl\{\eta(u)\varphi_t + q(u)\varphi_x\bigr\}\,dx\,dt \\[4mm]
~\ds=~-\int\sum_{\tau< t_m\leq\tau'} \Big(\eta(u (t_{m},x))-\eta(u(t_{m}-,x))\Big)\,\varphi(t_{m},x) dx\\[4mm]
~\ds\geq~-\int\sum_{\tau< t_m\leq\tau'} ~D\eta (u (t_{m},x))\big(u (t_{m},x)-u(t_{m}-,x)\big)\,\varphi(t_{m},x) dx\\[4mm]
~\ds\geq ~ -C_1\,\delta^2\,\Big({\tau'-\tau\over \ve}\Big)\Vert \varphi \Vert_{W^{1,\infty}}\, \hbox{Lip}(\eta)\cdot\Big(\sup_{t\in[\tau,\tau']} \hbox{Tot.Var.}\lbrace u(t, \cdot)\rbrace \Big)\\[4mm]
~\ds= ~ -C_2\,\ve\,({\tau'-\tau})\Vert \varphi \Vert_{W^{1,\infty}}\, \hbox{Lip}(\eta)\cdot\Big(\sup_{t\in[\tau,\tau']} \hbox{Tot.Var.}\lbrace u(t, \cdot)\rbrace \Big)\, ,
\enda$$
where the first inequality follows from the strict convexity of $\eta$,  by (\ref{ucx}).

\section{Numerical implementation}
\label{s:7}
\setcounter{equation}{0}
In this last section we discuss details of the post-processing algorithm, 
and present a numerical simulation.

STEP 1 of the algorithm, computing the total variation of the numerical solution $u(t,\cdot)$, is entirely straightforward.

STEP 2, identifying the location of the large shocks, requires more attention.
Given a pair of constants $K>0$ and $\sigma>\!>\ve>0$, we first identify regions  where
the total variation of $u$ is large.   For this purpose, we introduce
\begin{definition}\label{d:flag}
For a given function $u:[0,T]\times\R\mapsto\R^n$, the points $(t,x)$
such that
\bel{flag}
\min\Big\{ \tv\bigl\{u(t, \cdot)\,;~[x-\sigma, x+\ve]\bigr\}~,~~
 \tv\bigl\{u(t, \cdot)\,;~[x-\ve,\,x+\sigma]\bigr\}\Big\}~>~K\sigma\eeq
will be called {\bf flagged points}.
%Points where (\ref{flag}) does not hold will be called {\bf regular points}.
\end{definition}

Notice that, by definition, the 
oscillation of $u(t_j,\cdot)$ on a small interval 
to the left or to the right of a flagged point  must be large. 
Roughly speaking, the following result shows that, outside flagged points,  solutions are approximately Lipschitz continuous with constant $2K$.

\begin{lemma}\label{l:41}
Assume that all points $(t,x)$ with $x\in [a,b]$ are not flagged.
Then
\bel{nf}
|u(t,a)- u(t,b)|~\leq~ \left(1+{b-a\over \sigma+\ve}\right) 2\sigma K\,.\qquad \eeq
\end{lemma}

{\bf Proof.} {\bf 1.} By assumption, we can split the interval as 
$ [a,b]= I^+\cup I^-$, where
$I^+,I^-$ are two disjoint sets with the following property.
Setting
$$J_x~=~\left\{ \bega{rl} [x-\ve, x+\sigma] \quad &\hbox{if}\quad x\in I^+,\\[3mm]
 [x-\sigma, x+\ve] \quad &\hbox{if}\quad x\in I^-,\enda\right.$$
 one has
 $$\tv\{u(t,\cdot);~J_x\}~\leq~\sigma K\qquad \forall x\in [a,b].$$
\v
{\bf 2.} We claim that  every subinterval $[c,d]\subseteq [a,b]$
with length $d-c\leq \sigma+\ve$ can be covered by two of the intervals $J_x$.
Indeed, three cases can arise:

CASE 1: $c+\ve\in I^+$.   Then $[c,d]\subseteq J_{c+\ve}\,$.

CASE 2: $d-\ve\in I^-$.   Then $[c,d]\subseteq J_{d-\ve}\,$.

CASE 3: $c+\ve\in I^-$ and $d-\ve\in I^+$.   Then we can find
to points $c+\ve \leq x < y\leq d-\ve $, such that 
$$y-x< \ve,\qquad x\in I^-,~y\in I^+.$$
In this case, $[c,d]\subseteq J_x\cup J_y$, proving our claim.
\v
{\bf 3.} To complete the proof,
 we cover $[a,b]$ with finitely many  intervals $[c_i, d_i]$, $i=1,\ldots,N$, such that 
$$N~\leq~1+{b-a\over \sigma+\ve}\,,
\qquad \quad d_i-c_i ~\leq~\sigma+\ve\qquad\forall i=1,\ldots,N,$$
By the previous construction, for every $i$  we have
$$\tv\bigl\{u(t,\cdot)\,;~[c_i, d_i]\bigr\}~\leq~2K\sigma.$$
Therefore, the total variation of $u(t,\cdot)$ over $[a,b]$ is 
bounded by $2N\sigma K$. This yields (\ref{nf}).
\endproof

Having defined the set $\F\subseteq [0,T]\times\R$ of all flagged points,
for each $t_j= j\ve\in [0,T]$, $j=0,1,\ldots,\nu$, we denote by
$$\F_j~\doteq~\bigl\{ x\in \R\,;~~(t_j, x)\in \F\bigr\}$$
the set of flagged points at time $t_j$.

For every time $t_j$, we identify intervals $[a,b]$ such that  $b-a\leq \delta\doteq 
 \ve^{2/3}$ and moreover
$$ a\in \F_j\,,\qquad b\in \F_j\,,    \qquad \F_j\cap [a-\rho, a[~=~\F_j\,\cap\, ]b, b+\rho]~=~\emptyset.$$
In other words, the points $(t_j, a)$ and $(t_j,b)$ are flagged, but points to the left of $a$ and to the
right of $b$ are not flagged.   

Each such interval $[a,b]$ locates a possible isolated shock at time $t_j$.
To check if this shock can be traced over the entire interval $[t_j, t_{j+1}]$,
we check if there exists an interval $[c,d]$ with the same properties at time $t=t_{j+1}$,
with $$[c,d]~\subset~[a+\lambda^-h~, ~b+\lambda ^+ h].$$
In the positive case, we approximate the shock location as
$$\gamma(t)~=~x_0 + \lambda(t-t_j)$$
choosing $x_0, \lambda$ so that 
$$\gamma(t_j)~=~{a+b\over 2}\,,\qquad \qquad \gamma(t_{j+1})~=~{c+d\over 2}.$$

We then consider the polygonal regions $\Gamma$, $\Delta_l$, $\Delta'_r$ defined as in (\ref{GJI}).
If the two inequalities (\ref{osmall})-(\ref{jbg}) are both satisfied, we say that the 
parallelogram
$$\Gamma~=~\Big\{ (t,x)\,;~t\in [t_j, t_{j+1}]\,,~~x_0 -\delta+(t-t_j) \lambda\,\leq\,x\,\leq\, x_0+\delta +\lambda^-( t-t_j) \Big\},$$
traces the shock. 
The trapezoid
$$\Delta'~=~\Big\{ (t,x)\,;~t\in [t_j, t_{j+1}]\,,~~a'+\lambda^+(t-t_j)\,\leq\,x\,\leq\, b' +\lambda^-( t-t_j) \Big\},$$
$$a'~\doteq~x_0-\rho-\delta - (\lambda^+-\lambda^-)h,\qquad b'~\doteq~x_0+\rho+\delta + (\lambda^+-\lambda^-)h,$$
is then inserted within the list of trapezoids $\Delta^{(j\ell)}$ in  (\ref{Djl}), containing a traced shock.
On the other hand, if one of the inequalities (\ref{osmall})-(\ref{jbg}) fails, the shock is not traced.
\v
STEP 3 of the algorithm provides a covering of each domain  
$$\Big([t_j,t_{j+1}]\times \R\Big) \setminus \bigcup_ {\ell=1}^{N'(j)} \Delta^{(j\ell)},\qquad\quad j=0,1,\ldots, \nu-1,$$
with finitely many trapezoids $\Delta_{jk}$ as in (\ref{Tjk}).
This step is straightforward. The algorithm terminates
by computing the constants  $\kappa_j$ in (\ref{kjmax}), which  provide
an upper bound on the oscillation of $u$ on each $\Delta_{jk}$.
\v

\begin{figure}[ht]
\centerline{\hbox{\includegraphics[width=7cm]{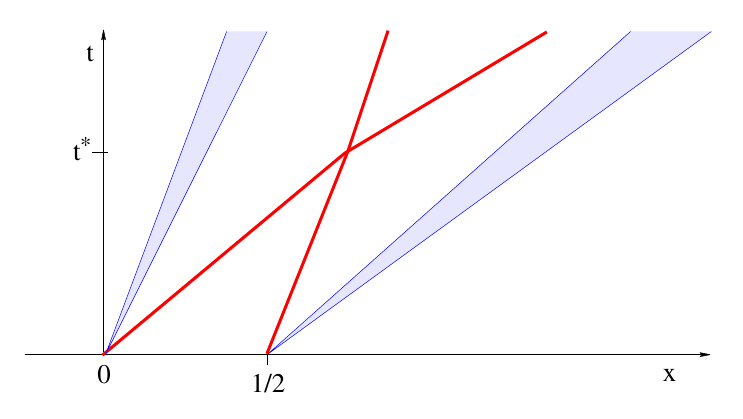}}}
\caption{\small A sketch of the exact solution to (\ref{GD})-(\ref{RD}), containing
two centered rarefaction waves and  two shocks, interacting at time $t^*$.}
\label{f:hyp225}
\end{figure}

{\bf Example.}
We consider  a model of
isentropic gas dynamics in Lagrangian coordinates.
Using a shifted system of coordinates, this can be written as
\bel{GD}
\left\{ \bega{rl} v_t - u_x +v_x&=~0,\\[3mm]
\ds u_t+\left({1\over 2 v^2}\right)_x+u_x&=0.\enda\right.\eeq
Here $u$ is the velocity of the gas, while $v$ denotes the specific volume.
By the choice of coordinates,
the characteristic speeds are
$$\lambda^\pm~=~1\pm v^{-3/2}.$$
In particular, when $v\geq  1$, one has $\lambda(v)\in[\lambda^-, \lambda^+]\doteq [0,2]$.
%\bel{lpmb}
%\lambda(v)~\in~[\lambda^-, \lambda^+]~\doteq~[0,2].\eeq
We consider the Cauchy problem with piecewise constant initial data
\bel{RD} v(0,x)~=~\left\{ \bega{cl} 2\quad&\hbox{if}~~x<0,\cr
3\quad &\hbox{if}~~0<x<{1/2},\cr
1\quad &\hbox{if}~~x>{1/2},\enda\right.\qquad\qquad u(0,x)~=~0.\eeq
The exact solution is shown in Fig.~\ref{f:hyp225}.

\begin{figure}[ht]
\centerline{\hbox{\includegraphics[width=8cm]{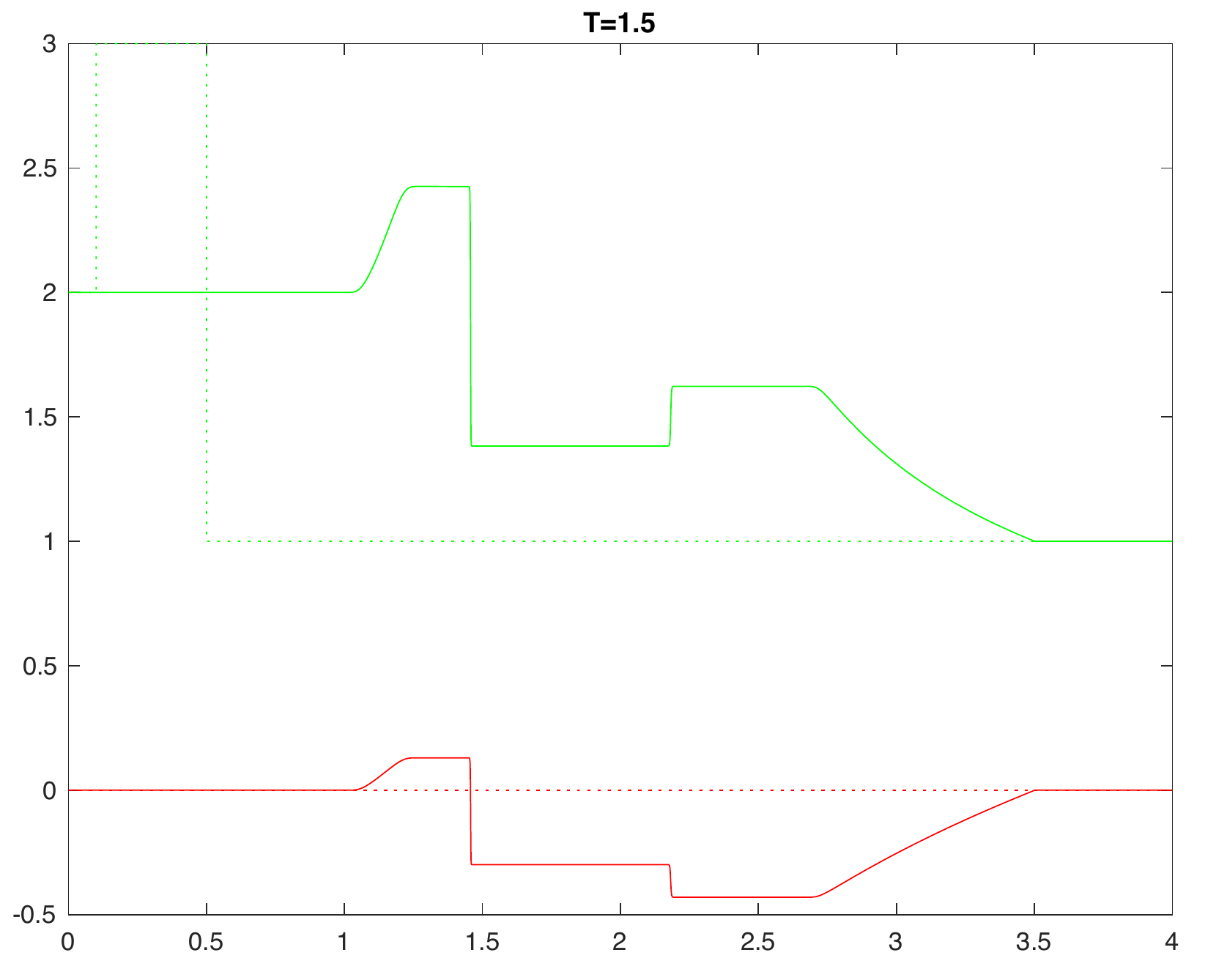}}}
\caption{\small The components of the solution at the terminal time $T=1.5$, computed by the Godunov
scheme.  Above: the specific volume $v(T,\cdot)$. Below: the velocity $u(T, \cdot)$.
}
\label{f:p-systUV1}
\end{figure}

We compute an approximate solution using
 the Godunov upwind scheme with mesh sizes
 $$\Delta x\,=\,\ve\,=\,0.0005,\qquad \qquad \Delta t\,=\,{\ve\over 2}\,=\,0.00025\,.$$
 The  profiles of the two components of the solution, at the terminal time $T=1.5$, are
 shown in Fig.~\ref{f:p-systUV1}.
 
\begin{figure}[ht]
\centerline{\hbox{\includegraphics[width=8cm]{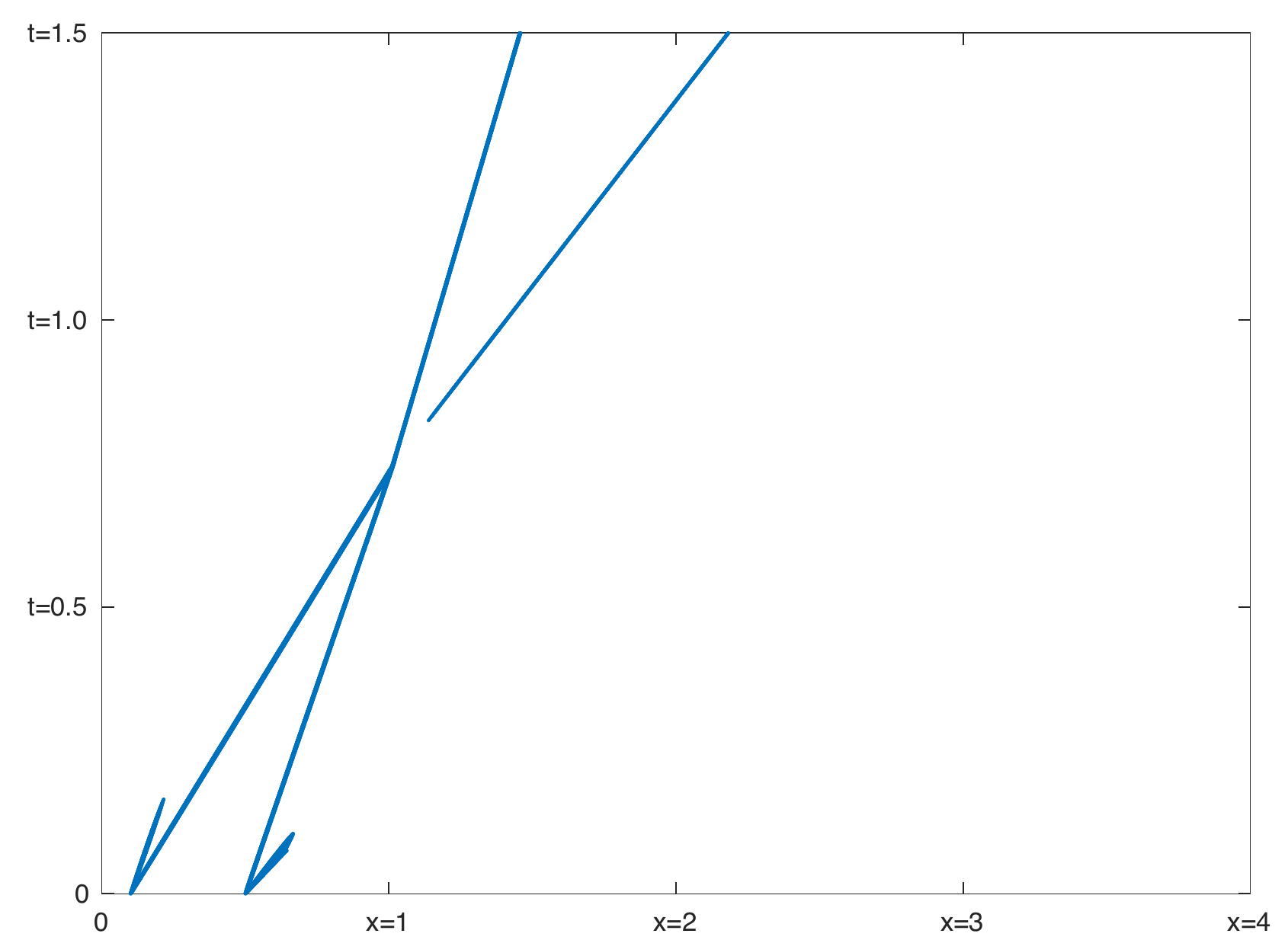}}\quad 
\hbox{\includegraphics[width=7.4cm]{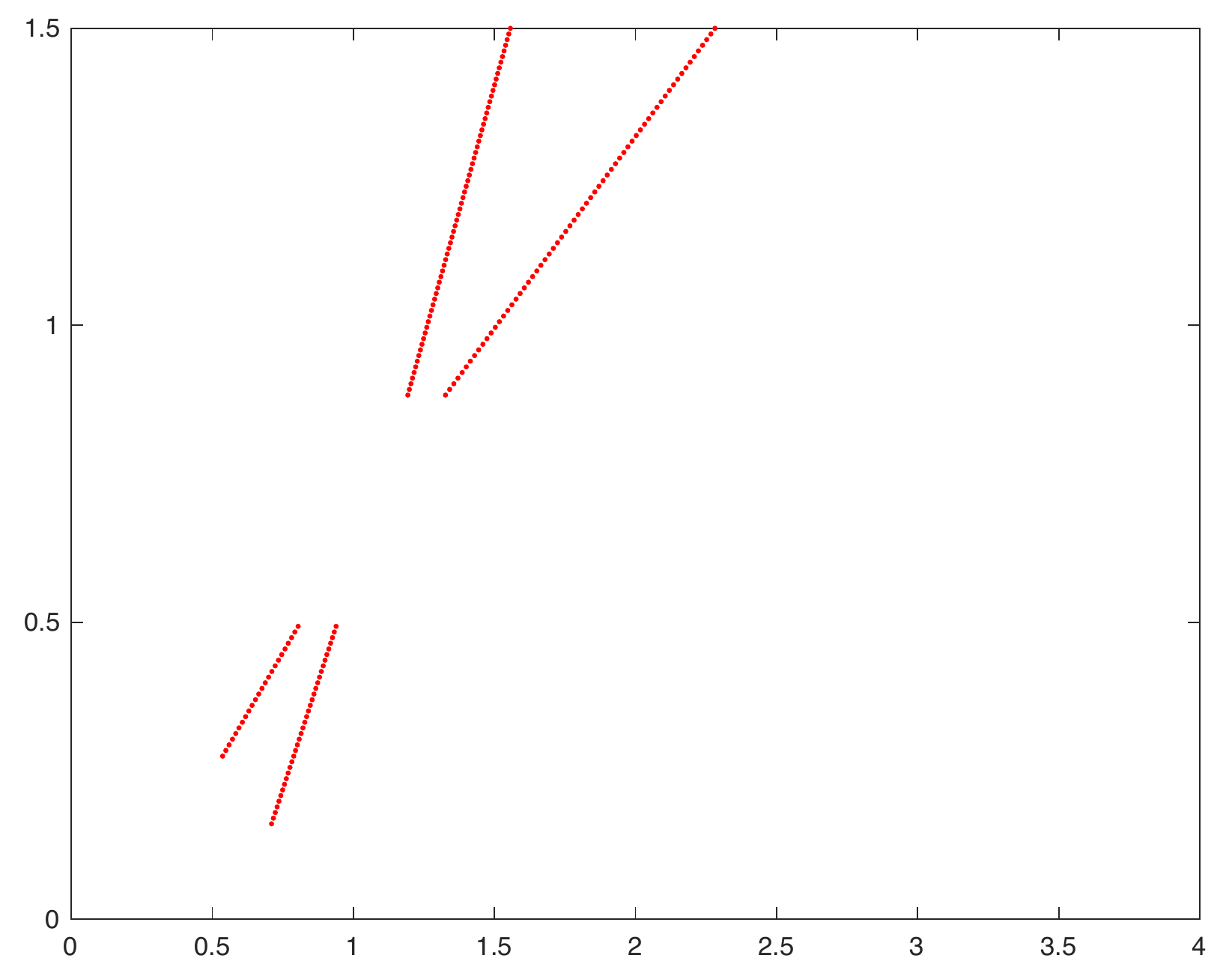}}}
\caption{\small Left: the points flagged by the post-processing algorithm.
Right: the portions of the two shocks that are actually traced.}
\label{f:flagpt1}
\end{figure}

To illustrate how the post-processing algorithm works, in Fig.~\ref{f:flagpt1}, left, 
we plot the set of flagged points.
These are computed according to Definition~\ref{d:flag},
choosing $\sigma= 0.0063$ and $K = 25$.
In Fig.~\ref{f:flagpt1}, right, we  identifying the shocks that can be traced on each 
time interval $[t_j, t_{j+1}]$. Here $\kappa' = 0.1$, $\sigma_{min}= 0.4$. 
Notice that, according to our previous construction, 
each trapezoid $\Delta^{(j\ell)}$ around a traced shock will have the form
$$\Delta^{(j\ell)}~=~\Big\{ (t,x)\,;~~t\in [t_j, \,t_{j+1}], ~~x_{j\ell} -\delta-\rho -2h + 2 (t-t_j) 
~\leq~ x~\leq~x_{j\ell} +\delta + \rho+ 2h\Big\}.$$
Indeed, this is obtained from (\ref{aabb})-(\ref{Dp}), with
$$a'\,=~\,x_{j\ell} -\delta-\rho -2h,\qquad b'\,=\,x_{j\ell} +\delta + \rho+ 2h,
\qquad \lambda^-=0,\qquad \lambda^+ =2.$$

\begin{figure}[ht]
\centerline{\hbox{\includegraphics[width=8cm]{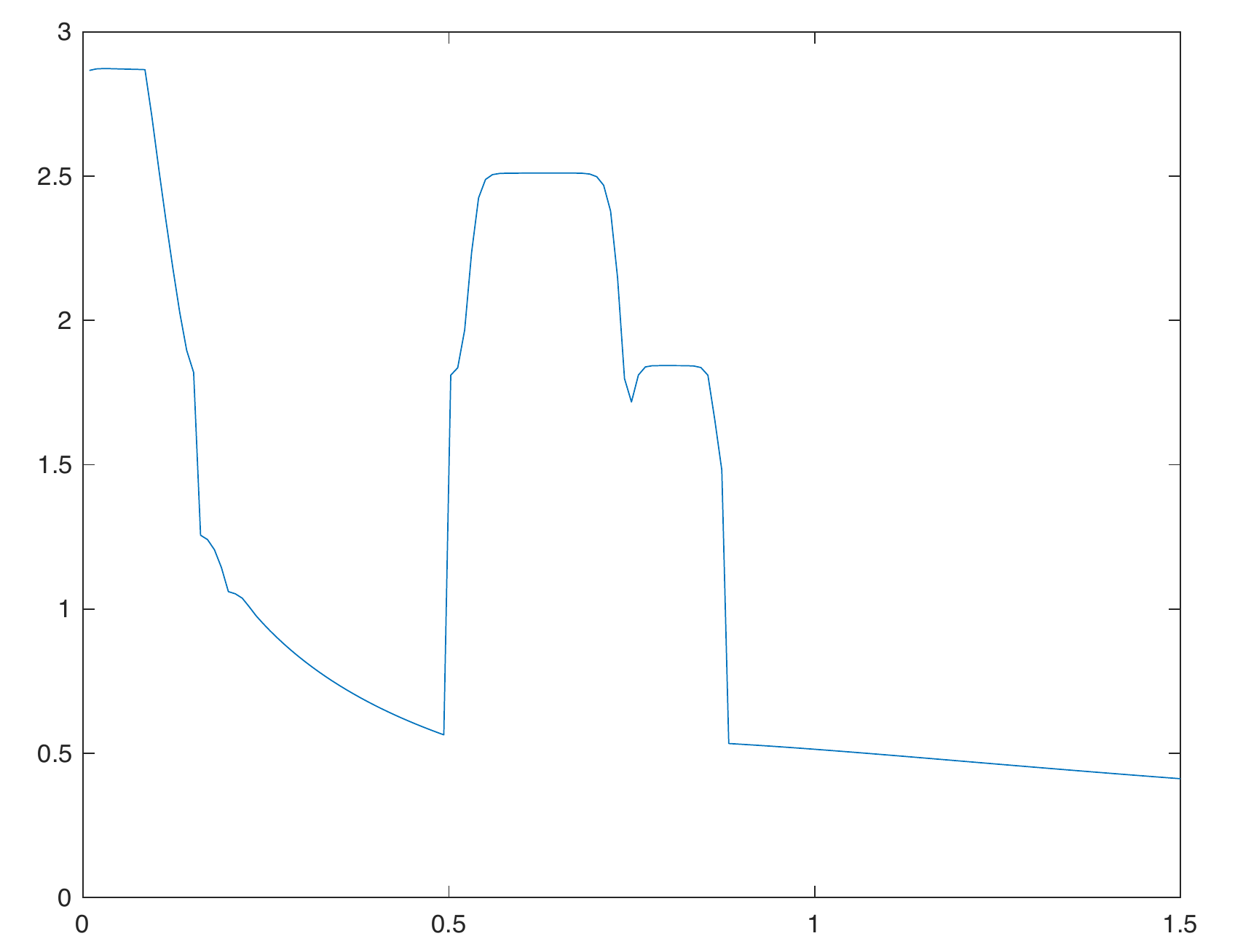}}}
\caption{\small An approximate computation of the function $\kappa(t)$ at (\ref{kap}),
determining the error rate.
}
\label{f:plotKK.pdf}
\end{figure}

Finally, in Fig.~\ref{f:plotKK.pdf} we plot an approximate graph of the function
\bel{kap}\kappa(t)~=~ \kappa_j~\doteq
~\max_{1\leq k\leq N(j)} \osc\{ u\,;~\Delta_{jk}\}\qquad \hbox{if}~~ t\in [t_j, t_{j+1}\,[\,.\eeq
One can think of $\kappa(t)$ as the maximum oscillation of the numerical solution
$u(t,\cdot)$ on domains of diameter $\O(1)\cdot  \ve^{1/3}$, outside the large traced shocks.
In view of (\ref{EEE}), this function $\kappa(\cdot)$ determines the  rate at which the 
distance $\|u(t,\cdot) - S_t \bar u\|_{\L^1}$ between the approximate and the exact solution
increases in time.
Notice that $\kappa(t)$ is large for $t\approx 0$, when the main contribution to the error
comes from the two centered rarefactions.  As time increases, the rarefactions decay, and 
the value of $\kappa(t)$ decays as well.   As $t$ approaches the interaction time $t^*$,
the two shocks cannot be individually traced. As a consequence, 
the  value of $\kappa(t)$ suddenly becomes very large.  
Finally, when the shocks 
move away from each other and can  be  traced once again, we see that 
$\kappa(t)$ regains its small values.

\v
{\bf Acknowledgments.}
This research  was partially supported by NSF with  
grant  DMS-2006884, ``Singularities and error bounds for hyperbolic equations".
\v

\end{document}